\documentclass[10pt]{amsart}
\input{diagrams}

\newtheorem{proposition}{Proposition}[section]
\newtheorem{lemma}[proposition]{Lemma}
\newtheorem{corollary}[proposition]{Corollary}
\newtheorem{theorem}[proposition]{Theorem}

\theoremstyle{definition}
\newtheorem{definition}[proposition]{Definition}
\newtheorem{example}[proposition]{Example}

\theoremstyle{remark}

\newcommand{\thlabel}[1]{\label{th:#1}}
\newcommand{\thref}[1]{Theorem~\ref{th:#1}}
\newcommand{\selabel}[1]{\label{se:#1}}
\newcommand{\seref}[1]{Section~\ref{se:#1}}
\newcommand{\lelabel}[1]{\label{le:#1}}
\newcommand{\leref}[1]{Lemma~\ref{le:#1}}
\newcommand{\prlabel}[1]{\label{pr:#1}}
\newcommand{\prref}[1]{Proposition~\ref{pr:#1}}
\newcommand{\colabel}[1]{\label{co:#1}}

\newcommand{\exlabel}[1]{\label{ex:#1}}

\newcommand{\delabel}[1]{\label{de:#1}}

\newcommand{\eqlabel}[1]{\label{eq:#1}}
\newcommand{\equref}[1]{(\ref{eq:#1})}

\def\equal#1{\smash{\mathop{=}\limits^{#1}}}

\newcommand{\Hom}{{\rm Hom}}

\newcommand{\Ker}{{\rm Ker}\,}

\newcommand{\im}{{\rm Im}\,}

\def\lan{\langle}
\def\ran{\rangle}
\def\ot{\otimes}

\newcommand{\Cc}{\mathcal{C}}

\newcommand{\Mm}{\mathcal{M}}

\newcommand{\Zz}{\mathcal{Z}}
\newcommand{\Ww}{\mathcal{W}}
\newcommand{\YD}{\mathcal{YD}}

\def\leftact{\hbox{$\rightharpoonup$}}
\def\rightact{\hbox{$\leftharpoonup$}}

\def\text#1{{\rm {\rm #1}}}

\def\ol{\overline}
\def\ul{\underline}

\begin{document}
\title[Yetter-Drinfeld modules over weak bialgebras]{Yetter-Drinfeld modules over weak bialgebras}
\author{S. Caenepeel}
\address{Faculty of Engineering Sciences,
Vrije Universiteit Brussel, VUB, B-1050 Brussels, Belgium}
\email{scaenepe@vub.ac.be}
\urladdr{http://homepages.vub.ac.be/\~{}scaenepe/}
\author{Dingguo Wang}
\address{Department of Mathematics, Qufu Normal University,
Qufu, Shandong 273165, China}
\email{dgwang@qfnu.edu.cn}
\author{Yanmin Yin}
\address{Department of Mathematics, Shandong Institute of Architecture and
Engineering, Jinan, Shandong 250014, China}
\email{yanmin\hbox{\ul{\hspace*{2mm}}}yin@163.com}
\thanks{Research supported by the projects G.0278.01 ``Construction
and applications of non-commutative geometry: from algebra to physics"
from FWO-Vlaanderen and ``New computational, geometric and algebraic methods applied to quantum groups and differential operators" from the Flemish and Chinese
governments.}

\subjclass{16W30}
\keywords{}
\begin{abstract}
We discuss properties of Yetter-Drinfeld modules over weak bialgebras
over commutative rings. The categories of left-left, left-right, right-left
and right-right Yetter-Drinfeld modules over a weak Hopf algebra are isomorphic as braided monoidal
categories. Yetter-Drinfeld modules can be viewed as weak Doi-Hopf modules,
and, a fortiori, as weak entwined modules. If $H$ is finitely generated and
projective, then we introduce the Drinfeld double using duality results between
entwining structures and smash product structures, and show that the category
of Yetter-Drinfeld modules is isomorphic to the category of modules over the
Drinfeld double. The category of finitely generated projective Yetter-Drinfeld modules
over a weak Hopf algebra has duality.
\end{abstract}

\maketitle

\section*{Introduction}
Weak bialgebras and Hopf algebras are generalizations of ordinary bialgebras
and Hopf algebras in the following sense: the defining axioms are the same,
but the multiplicativity of the counit and comultiplicativity of the unit
are replaced by weaker axioms. The easiest example of a weak Hopf algebra
is a groupoid algebra; other examples are face algebras \cite{Hayashi},
quantum groupoids \cite{Ocneanu}, generalized Kac algebras \cite{Yamanouchi} and
quantum transformation groupoids \cite{NV}. Temperley-Lieb algebras give rise
to weak Hopf algebras (see \cite{NV}).
A purely algebraic study of weak Hopf algebras has been presented in
\cite{BohmNSI}. A survey of weak Hopf algebras and their applications
may be found in \cite{NV}. It has turned out that many results of
classical Hopf algebra theory can be generalized to weak Hopf algebras.\\
Yetter-Drinfeld modules over finite dimensional weak Hopf algebras over
fields have been introduced by Nenciu \cite{Nenciu}. It is shown in \cite{Nenciu}
that the category of finite dimensional Yetter-Drinfeld modules is
isomorphic to the category of finite dimensional modules over the Drinfeld
double, as introduced in the appendix of \cite{Bohm}. It is also shown that
this category is braided isomorphic to the center of the category of finite dimensional
$H$-modules. In this note, we discuss Yetter-Drinfeld modules over weak
bialgebras over commutative rings. The results in \cite{Nenciu} are slightly
generalized and more properties are given.\\
In \seref{2}, we compute the weak center of the category of modules over
a weak bialgebra $H$, and show that it is isomorphic to the category of Yetter-Drinfeld
modules. If $H$ is a weak Hopf algebra, then the weak center equals the center.
In this situation, properties of the center construction can be applied to show that
the four categories of Yetter-Drinfeld modules, namely
the left-left, left-right, right-left and right-right versions, are isomorphic
as braided monoidal categories. Here we apply methods that have been used
before in \cite{BCP}, in the case of quasi-Hopf algebras.\\
In \cite{cmz}, it was observed that Yetter-Drinfeld modules over a classical
Hopf algebra are special cases of Doi-Hopf modules, as introduced by
Doi and Koppinen (see \cite{Doi92,Koppinen95}). In \seref{3}, we will show
that Yetter-Drinfeld modules over weak Hopf algebras are weak Doi-Hopf modules,
in the sense of B\"ohm \cite{Bohm}, and, a fortiori, weak entwined modules
\cite{CDG}, and comodules over a coring \cite{BW}.\\
The advantage of this approach is that it leads easily to a new description
of the Drinfeld double of a finitely generated projective weak Hopf algebra,
using methods developed in \cite{CDG}: we define the Drinfeld double as
a weak smash product of $H$ and its dual. We show that our Drinfeld double is
equal to the Drinfeld double of \cite{Bohm,Nenciu} (see \prref{4.2}) and anti-isomorphic
to the Drinfeld double of \cite{NVT} (see \prref{4.4}). In \seref{5}, we show that the
category of finitely generated projective Yetter-Drinfeld modules over a weak
Hopf algebra has duality.\\
In Sections \ref{se:1.1} and \ref{se:1.2}, we recall some general properties
of weak bialgebras and Hopf algebras. Further detail can be found in
\cite{BW,BohmNSI,NV}.
In \seref{1.3}, we recall the center
construction, and in \seref{1.4}, we recall the notions of weak Doi-Hopf
modules, weak entwining structures and weak smash products.

\section{Preliminary results}\selabel{1}
\subsection{Weak bialgebras}\selabel{1.1}
Let $k$ be a commutative ring. Recall that a weak $k$-bialgebra is a $k$-module
with a $k$-algebra structure $(\mu,\eta)$ and a $k$-coalgebra structure
$(\Delta,\varepsilon)$ such that
$\Delta(hk)=\Delta(h)\Delta(k)$,
for all $h,k\in H$, and
\begin{eqnarray}
\Delta^2(1)&=& 1_{(1)}\ot 1_{(2)}1_{(1')}\ot 1_{(2')}=
1_{(1)}\ot 1_{(1')}1_{(2)}\ot 1_{(2')},\eqlabel{1.1.1}\\
\varepsilon(hkl)&=& \varepsilon(hk_{(1)})\varepsilon(k_{(2)}l)=
\varepsilon(hk_{(2)})\varepsilon(k_{(1)}l),\eqlabel{1.1.2}
\end{eqnarray}
for all $h,k,l\in H$. We use the Sweedler-Heyneman notation for the
comultiplication, namely
$$\Delta(h)=h_{(1)}\ot h_{(2)}=h_{(1')}\ot h_{(2')}.$$
We summarize the elementary properties of weak bialgebras. The proofs are
direct applications of the defining axioms (see \cite{BohmNSI,NV}).
We have idempotent maps $\varepsilon_t,~\varepsilon_s:\ H\to H$ defined by
$$\varepsilon_t(h)=\varepsilon(1_{(1)}h)1_{(2)}~~;
\varepsilon_s(h)=1_{(1)}\varepsilon(h1_{(2)}).$$
$\varepsilon_t$ and $\varepsilon_s$ are called the target map and the
source map, and their images $H_t=\im(\varepsilon_t)=\Ker(H-\varepsilon_t)$ and
$H_s=\im(\varepsilon_s)=\Ker(H-\varepsilon_s)$ are called the target and source space. For all
$g,h\in H$, we have
\begin{equation}\eqlabel{1.1.3}
h_{(1)}\ot \varepsilon_t(h_{(2)})=1_{(1)}h\ot 1_{(2)}~~{\rm and}~~
\varepsilon_s(h_{(1)})\ot h_{(2)}=1_{(1)}\ot h1_{(2)},
\end{equation}
and
\begin{equation}\eqlabel{1.1.3*}
h\varepsilon_t(g)=\varepsilon(h_{(1)}g)h_{(2)}~~{\rm and}~~
\varepsilon_s(g)h=h_{(1)}\varepsilon(gh_{(2)}).
\end{equation}
From \equref{1.1.3*}, it follows immediately that
\begin{equation}\eqlabel{1.1.3a}
\varepsilon(h\varepsilon_t(g))=\varepsilon(hg)~~{\rm and}~~
\varepsilon(\varepsilon_s(g)h)=\varepsilon(gh).
\end{equation}
The source and target space can be described as follows:
\begin{eqnarray}
H_t=\{h\in H~|~\Delta(h)=1_{(1)}h\ot 1_{(2)}\}=
\{\phi(1_{(1)})1_{(2)}~|~\phi\in H^*\};\eqlabel{1.1.4}\\
H_s=\{h\in H~|~\Delta(h)=1_{(1)}\ot h1_{(2)}\}=
\{1_{(1)}\phi(1_{(2)})~|~\phi\in H^*\}.\eqlabel{1.1.5}
\end{eqnarray}
We also have
\begin{equation}\eqlabel{1.1.6}
\varepsilon_t(h)\varepsilon_s(k)=\varepsilon_s(k)\varepsilon_t(h),
\end{equation}
and its dual property
\begin{equation}\eqlabel{1.1.6*}
\varepsilon_s(h_{(1)})\ot \varepsilon_t(h_{(2)})=
\varepsilon_s(h_{(2)})\ot \varepsilon_t(h_{(1)}).
\end{equation}
Finally $\varepsilon_s(1)=\varepsilon_t(1)=1$, and
\begin{equation}\eqlabel{1.1.7}
\varepsilon_t(h)\varepsilon_t(g)=\varepsilon_t(\varepsilon_t(h)g)~~{\rm and}~~
\varepsilon_s(h)\varepsilon_s(g)=\varepsilon_s(h\varepsilon_s(g)).
\end{equation}
This implies that $H_s$ and $H_t$ are subalgebras of $H$.\\

\begin{lemma}\lelabel{1.2}
Let $H$ be a weak bialgebra over a commutative ring. Then
$\Delta(1)\in H_s\ot H_t$.
\end{lemma}

\begin{proof}
Applying $H\ot \varepsilon\ot H$ to \equref{1.1.1},
we find that $1_{(1)}\ot 1_{(2)}=\varepsilon_s(1_{(1)})\ot 1_{(2)}
\in H_s\ot H$ and
$1_{(1)}\ot 1_{(2)}=1_{(1)}\ot \varepsilon_t(1_{(2)})\in H\ot H_t$.
Now let $K_s=\Ker(\varepsilon_s)$, $K_t=\Ker(\varepsilon_t)$. Then
$H= H_s\oplus K_s= H_t\oplus K_t,$
and
$$H\ot H=H_s\ot H_t\oplus H_s\ot K_t \oplus K_s\ot H_t \oplus K_s\ot K_t,$$
so it follows that $H_s\ot H_t= H\ot H_t \cap H_s\ot H$.
\end{proof}

The target and source map for the weak bialgebra $H^{\rm op}$ are
\begin{equation}\eqlabel{1.2a.1}
\ol{\varepsilon}_t(h)=\varepsilon(h1_{(1)})1_{(2)}\in H_t~~{\rm and}~~
\ol{\varepsilon}_s(h)=\varepsilon(1_{(2)}h)1_{(1)}\in H_s.
\end{equation}
$\ol{\varepsilon}_t$ and $\ol{\varepsilon}_s$ are also projections.\\

The source and target space are anti-isomorphic, and they are separable Frobenius
algebras over $k$. This was first proved for weak Hopf algebras (see \cite{BohmNSI}), 
and then generalized to weak bialgebras (see \cite{Schauenburg}).

\begin{lemma}\eqlabel{1.2a} \cite{Schauenburg}
Let $H$ be a weak bialgebra. Then $\ol{\varepsilon}_s$ restricts to an anti-algebra
isomorphism $H_t\to H_s$ with inverse $\varepsilon_t$, and 
$\ol{\varepsilon}_t$ restricts to an anti-algebra
isomorphism $H_s\to H_t$ with inverse $\varepsilon_s$.
\end{lemma}

\begin{proposition}\prlabel{1.2b} \cite{Schauenburg}
Let $H$ be a weak bialgebra. Then $H_s$ and $H_t$ are Frobenius separable
$k$-algebras. The separability idempotents of $H_t$ and $H_s$ are
$$e_t=\varepsilon_t(1_{(1)})\ot 1_{(2)}=1_{(2)}\ot \ol{\varepsilon}_t(1_{(1)});$$
$$e_s=1_{(1)}\ot \varepsilon_s(1_{(2)})=\ol{\varepsilon}_s(1_{(2)})\ot 1_{(1)}.$$
The Frobenius systems for $H_t$ and $H_s$ are respectively
$(e_t,\varepsilon_{|H_t})$ and $(e_s, \varepsilon_{|H_s})$. In particular, we have
for all $z\in H_t$ that
\begin{equation}\eqlabel{1.2b.1}
z\varepsilon_t(1_{(1)})\ot 1_{(2)}=
\varepsilon_t(1_{(1)})\ot 1_{(2)}z.
\end{equation}
\end{proposition}

It was shown in \cite{NVT} that 
the category of modules over a weak Hopf algebra is monoidal; it follows from
the results of \cite{Schauenburg} that this property can be generalized to weak
bialgebras. We explain now how this can be done directly.\\
Let $M$ be a left $H$-module. By restriction of scalars, $M$ is a left $H_t$-module;
$M$ becomes an $H_t$-bimodule, if we define a right $H_t$-action by
$$m\cdot z=\ol{\varepsilon}_s(z)m.$$\\
Let $M,N\in {}_H\Mm$, the category of left $H$-modules. We define
$$M\ot_t N=\Delta(1)(M\ot N),$$
the $k$-submodule of $M\ot N$ generated by elements of the form
$1_{(1)}\ot 1_{(2)}$. $M\ot_t N$ is a left $H$-module, with left diagonal
action $h\cdot(m\ot n)=h_{(1)}m\ot h_{(2)}n$. It follows from
\equref{1.1.1} that the tensor product $\ot_t$ is associative. Observe
that
$$M\ot_t N\ot_t P=\Delta^2(1)(M\ot N\ot P).$$
$H_t\in {}_H\Mm$, with
left $H$-action $h\leftact z=\varepsilon_t(hz)$. The induced $H_t$-bimodule
structure is given by left and right multiplication by elements of $H_t$.\\
For $M,N\in {}_H\Mm$, consider the projection
$$\pi:\ M\ot N\to M\ot_t N,~~\pi(m\ot n)=1_{(1)}m\ot 1_{(2)}n.$$
Applying $\ol{\varepsilon}_s\ot H_t$ to \equref{1.2b.1}, we find
$$\ol{\varepsilon}_s(z\varepsilon_t(1_{(1)}))\ot 1_{(2)}=
1_{(1)}\ol{\varepsilon}_s(z)\ot 1_{(2)}=1_{(1)}\ot 1_{(2)}z,$$
hence
$$
\pi(mz\ot n)=\pi(\ol{\varepsilon}_s(z)m\ot n)=
1_{(1)}\ol{\varepsilon}_s(z)m\ot 1_{(2)}n
=1_{(1)}m\ot 1_{(2)}zn=\pi(m\ot zn).$$
So $\pi$ induces a map $\ol{\pi}:\ M\ot_{H_t} N\to M\ot_t N$, which is a left $H_t$-module
isomorphism with inverse given by
$$\ol{\pi}^{-1}(1_{(1)}m\ot 1_{(2)}n)=1_{(1)}m\ot_{H_t} 1_{(2)}n=m\ot_{H_t}n.$$

\begin{proposition}\prlabel{1.2c}
Let $H$ be a weak bialgebra. Then we have a monoidal category
$({}_H\Mm,\ot_t,H_t,a,l,r)$. The associativity constraints are the natural ones. The left and right
unit constraints $l_M:\ H_t\ot_t M\to M$ and $r_M:\ M\ot_t H_t\to M$ and their inverses
are given by the formulas
$$l_M(1_{(1)}\leftact z\ot 1_{(2)}m)=zm~~;~~l_M^{-1}(m)=\varepsilon_t(1_{(1)})\ot 1_{(2)}m;$$
$$r_M(1_{(1)}m\ot 1_{(2)}\leftact z)=\ol{\varepsilon}_s(z)m~~;~~
r_M^{-1}(m)=1_{(1)}m\ot 1_{(2)}.$$
\end{proposition}

\begin{proof}
This is a direct consequence of the observations made above. Let us check that
\begin{eqnarray*}
&&\hspace*{-2cm}
l_M^{-1}(l_M(1_{(1)}\leftact z\ot 1_{(2)}m))=l_M^{-1}(zm)=
\varepsilon_t(1_{(1)})\ot 1_{(2)}zm\\
&=&z\varepsilon_t(1_{(1)})\ot 1_{(2)}m\equal{\equref{1.1.7}}\varepsilon_t(z1_{(1)})\ot 1_{(2)}m\\
&=&\varepsilon_t(1_{(1)}z)\ot 1_{(2)}m=1_{(1)}\leftact z\ot 1_{(2)}m\\
&&\hspace*{-2cm}
l_M(l_M^{-1}(m))=l_M(\varepsilon_t(1_{(1)})\ot 1_{(2)}m)=m\\
&&\hspace*{-2cm}
r_M^{-1}(r_M(1_{(1)}m\ot 1_{(2)}\leftact z))=r_M^{-1}(\ol{\varepsilon}_s(z)m)=
1_{(1)}\ol{\varepsilon}_s(z)m\ot 1_{(2)}\\
&=& 1_{(1)}m\ot 1_{(2)}z=1_{(1)}m\ot 1_{(2)}\leftact z\\
&&\hspace*{-2cm}
r_M(r_M^{-1}(m))=r_M(1_{(1)}m\ot 1_{(2)})=\ol{\varepsilon}_s(1)m=m.
\end{eqnarray*}
\end{proof}

\subsection{Weak Hopf algebras}\selabel{1.2}
A weak Hopf algebra is a weak bialgebra together with a map $S:\ H\to H$,
called the antipode, satisfying
\begin{equation}\eqlabel{1.3.1}
S*H=\varepsilon_s,~~H*S=\varepsilon_t,~~{\rm and}~~S*H*S=S,
\end{equation}
where $*$ is the convolution product. It follows immediately that
\begin{equation}\eqlabel{1.3.2}
S=\varepsilon_s*S=S*\varepsilon_t.
\end{equation}
If the antipode exists, then it is unique. We will always
assume that $S$ is bijective; if $H$ is a finite dimensional weak Hopf algebra
over a field, then $S$ is automatically bijective (see \cite[Theorem 2.10]{BohmNSI}).

\begin{lemma}\lelabel{1.5}
Let $H$ be a weak Hopf algebra. Then $S$ is an anti-algebra and an anti-coalgebra
morphism.
For all $h,g\in H$, we have
\begin{eqnarray}
\varepsilon_t(hg)&=&
\varepsilon_t(h\varepsilon_t(g))=h_{(1)}\varepsilon_t(g)S(h_{(2)});\eqlabel{1.4.1}\\
\varepsilon_s(hg)&=&
\varepsilon_s(\varepsilon_s(h)g)=S(g_{(1)})\varepsilon_s(h)g_{(2)};
\eqlabel{1.4.1a}\\
\Delta(\varepsilon_t(h))&=&h_{(1)}S(h_{(3)})\ot\varepsilon_t(h_{(2)})\eqlabel{1.4.2}\\
\Delta(\varepsilon_s(h))&=&\varepsilon_s(h_{(2)})\ot S(h_{(1)})h_{(3)}\eqlabel{1.4.2a}.
\end{eqnarray}
\end{lemma}

\begin{lemma}\lelabel{1.6}
Let $H$ be a weak Hopf algebra. For all $h\in H$, we have
\begin{eqnarray}
\varepsilon_t(h)&=&\varepsilon(S(h)1_{(1)})1_{(2)}=
\varepsilon(1_{(2)}h)S(1_{(1)})=S(\ol{\varepsilon}_s(h))\eqlabel{1.6.1}\\
\varepsilon_s(h)&=&1_{(1)}\varepsilon(1_{(2)}S(h))=
\varepsilon(h1_{(1)})S(1_{(2)})=S(\ol{\varepsilon}_t(h))\eqlabel{1.6.2}.
\end{eqnarray}
\end{lemma}

\begin{corollary}\colabel{1.6b}
Let $H$ be a weak Hopf algebra. For all $h\in H$, we have
\begin{equation}\eqlabel{1.6.3}
\varepsilon_t(h_{(1)})\ot h_{(2)}=S(1_{(1)})\ot 1_{(2)}h~~;~~
h_{(1)}\ot \varepsilon_s(h_{(2)})=h1_{(1)}\ot S(1_{(2)}).
\end{equation}
\end{corollary}

\begin{proposition}\prlabel{1.7}
Let $H$ be a weak Hopf algebra. Then
\begin{equation}\eqlabel{1.7.1}
\varepsilon_t\circ S=\varepsilon_t\circ \varepsilon_s= S\circ \varepsilon_s~~;~~
\varepsilon_s\circ S=\varepsilon_s\circ \varepsilon_t= S\circ \varepsilon_t.
\end{equation}
\end{proposition}

\begin{corollary}\colabel{1.8}
Let $H$ be a weak Hopf algebra with bijective antipode. Then 
$S_{|H_t}=({\varepsilon}_s)_{|H_t}$, and $S^{-1}_{|H_s}=(\ol{\varepsilon}_t)_{|H_s}$,
so $S$ restricts
to an anti-algebra isomorphism $H_t\to H_s$.
\end{corollary}

It follows that the separability idempotents of $H_t$ and $H_s$ are
$e_t=S(1_{(1)})\ot 1_{(2)}$ and $e_s=1_{(1)}\ot S(1_{(2)})$. Consequently, we
have the following formulas, for $z\in H_t$ and $y\in H_s$:
\begin{eqnarray}
zS(1_{(1)})\ot 1_{(2)}&=&S(1_{(1)})\ot 1_{(2)}z;\eqlabel{1.9.3}\\
y1_{(1)}\ot 1_{(2)}&=& 1_{(1)}\ot S^{-1}(y)1_{(2)}.\eqlabel{1.9.4}
\end{eqnarray}
Applying $S^{-1}\ot H$ to \equref{1.9.3}, we find
\begin{equation}\eqlabel{1.9.2}
1_{(1)}S^{-1}(z)\ot 1_{(2)}=1_{(1)}\ot 1_{(2)}z.
\end{equation}

\subsection{The center of a monoidal category}\selabel{1.3}
Let $\Cc=(\Cc,\ot,I,a,l,r)$ be a monoidal category. 
The weak left center $\Ww_l(\Cc)$ is the category
with the following objects and morphisms. An object is a couple $(M,\sigma_{M, -})$,
with $M\in \Cc$ and $\sigma_{M, -}:\ M\ot -\to -\ot M$ a natural transformation, 
satisfying the following condition, for all $X,Y\in \Cc$:
\begin{equation}\eqlabel{center1}
(X\ot \sigma_{M, Y})\circ a_{X,M,Y}\circ (\sigma_{M, X}\ot Y)=a_{X,Y,M}\circ \sigma_{M, X\ot Y}
\circ a_{M,X,Y},
\end{equation}
and such that $\sigma_{M, I}$ is the composition of the natural isomorphisms
$M\ot I\cong M\cong I\ot M$. A morphism between $(M,\sigma_{M, -})$ and $(M',\sigma_{M', -})$ 
consists of $\vartheta:\ M\to M'$ in $\Cc$ such that
$$
(X\ot \vartheta)\circ \sigma_{M, X}=\sigma_{M', X}\circ (\vartheta\ot X).
$$
The left center $\Zz_l(\Cc)$ is the full subcategory of $\Ww_l(\Cc)$ consisting of
objects $(M,\sigma_{M, -})$ with $\sigma_{M,-}$ a natural isomorphism.
$\Zz_l(\Cc)$ is a braided monoidal category. 
The tensor product is
$$
(M,\sigma_{M, -})\ot (M',\sigma_{M', -})=(M\ot M', \sigma_{M\ot M', -})
$$
with
\begin{equation}\eqlabel{center2}
\sigma_{M\ot M', X}=a_{X,M,M'}\circ (\sigma_{M, X}\ot M')\circ a^{-1}_{M,X,M'}
\circ (M\ot \sigma_{M', X})\circ a_{M,M',X},
\end{equation}
and the unit is $(I,\sigma_{I,-})$, with
\begin{equation}\eqlabel{center2bis}
\sigma_{I,M}=r_M^{-1}\circ l_M.
\end{equation}
The braiding $c$ on $\Zz_l(\Cc)$ is given by
\begin{equation}\eqlabel{center3}
c_{M,M'}=\sigma_{M, M'}:\ (M, \sigma_{M, -})\ot (M', \sigma_{M', -})
\to (M',\sigma_{M', -})\ot (M,\sigma_{M, -}).
\end{equation}
$\Zz_l(\Cc)^{\rm in}$ will be our notation for the
monoidal category $\Zz_l(\Cc)$, together with the inverse braiding
$\tilde{c}$ given by $\tilde{c}_{M,M'}=c^{-1}_{M',M}=\sigma^{-1}_{M', M}$.

The right center $\Zz_r(\Cc)$ is defined in a similar way. An object is a couple 
$(M,\tau_{-, M})$, where $M\in \Cc$ and
$\tau_{-, M}:\ -\ot M\to M\ot -$ is a family of natural 
isomorphisms such that $\tau_{-, I}$ is the natural isomorphism 
and 
\begin{equation}\eqlabel{center4}
a^{-1}_{M, X, Y}\circ \tau_{X\ot Y, M}\circ a^{-1}_{X, Y, M}=(\tau_{X, M}\ot Y)
\circ a^{-1}_{X, M, Y}\circ (X\ot \tau_{Y, M}),
\end{equation}
for all $X,Y\in \Cc$. 
A morphism between $(M,\tau_{-, M})$ and $(M',\tau_{-, M'})$ consists of 
$\vartheta:\ M\to M'$ in $\Cc$ such that
$$
(\vartheta\ot X)\circ \tau_{X, M}=\tau_{X, M'}\circ (X\ot \vartheta),
$$
for all $X\in \Cc$.  
$\Zz_r(\Cc)$ is a braided monoidal category. The unit is $(I, l_{-}^{-1}\circ r_{-})$ and the 
tensor product is
$$
(M,\tau_{-, M})\ot (M',\tau_{-, M'})=(M\ot M', \tau_{-, M\ot M'})
$$
with
\begin{equation}\eqlabel{center5}
\tau_{X, M\ot M'}=a^{-1}_{M,M',X}\circ (M\ot \tau_{X, M'})\circ a_{M,X,M'}
\circ (\tau_{X, M}\ot M')\circ a^{-1}_{X,M,M'}.
\end{equation}
The braiding $d$ is given by
\begin{equation}\eqlabel{center6}
d_{M,M'}=\tau_{M, M'}:\ (M, \tau_{-, M})\ot (M', \tau_{-, M'})\to 
(M', \tau_{-, M'})\ot (M, \tau_{-, M}).
\end{equation}
$\Zz_r(\Cc)^{\rm in}$ is the monoidal category $\Zz_r(\Cc)$ with the inverse
braiding $\tilde{d}$ given by $\tilde{d}_{M,M'}=d^{-1}_{M',M}=\tau^{-1}_{M', M}$.

For details in the case where $\Cc$ is a strict monoidal category, we refer
to \cite[Theorem XIII.4.2]{Kassel}. The results remain valid in the case of
an arbitrary monoidal category, since every monoidal category is equivalent
to a strict one. Recall the following result from \cite{BCP}.

\begin{proposition}\prlabel{1.0.1}
Let $\Cc$ be a monoidal category. Then we have an isomorphism of
braided monoidal categories $F:\ \Zz_l(\Cc)\to \Zz_r(\Cc)^{\rm in}$, given by
$$
F(M, \sigma_{M, -})=(M, \sigma^{-1}_{M, -})~~{\rm and}~~F(\vartheta)=\vartheta .
$$
\end{proposition}

We have
a second monoidal structure on $\Cc$, defined as follows:
$$
\ol{\Cc}=(\Cc,\ol{\ot}=\ot\circ\tau,I, \ol{a}, r, l)
$$
with $\tau:\ \Cc\times\Cc\to \Cc\times \Cc$, $\tau(M,N)=(N,M)$ and $\ol{a}$ 
defined by $\ol{a}_{M,N,X}=a^{-1}_{X,N,M}$. 

If $c$ is a braiding on $\Cc$, then $\ol{c}$, given by $\ol{c}_{M,N}=
c_{N,M}$ is a braiding on $\ol{\Cc}$. In \cite{BCP}, the following
obvious result was stated.

\begin{proposition}\prlabel{1.0.2}
Let $\Cc$ be a monoidal category. Then 
$$
\ol{\Zz_l(\Cc)}\cong \Zz_r(\ol{\Cc})~~;~~
\ol{\Zz_r(\Cc)}\cong \Zz_l(\ol{\Cc})
$$ 
as braided monoidal categories.
\end{proposition}

\subsection{Weak entwining structures and weak smash products}\selabel{1.4}
The results in this Section are taken from \cite{CDG}. Let $A$ be a ring
without unit. $e\in A$ is called a preunit if $ea=ae=ae^2$, for all $a\in A$.
Then map $p:\ A\to A$, $p(a)=ae$, satisfies the following properties:
$p\circ p=p$ and $p(ab)=p(a)p(b)$. Then $\ol{A}={\rm Coim}(p)$ is a ring
with unit $\ol{e}$ and $\ul{A}=\im(p)$ is a ring with unit $e^2$. $p$ induces a
ring isomorphism $\ol{A}\to \ul{A}$.\\
Let $k$ be a commutative ring, $A$, $B$ $k$-algebras with unit, and
$R:\ B\ot A\to A\ot B$ a $k$-linear map. We use the notation
\begin{equation}\eqlabel{smash}
R(b\ot a)=a_R\ot b_R=a_r\ot b_r,
\end{equation}
where the summation is implicitely understood. $A\#_RB$ is the $k$-algebra
$A\ot B$ with newly defined multiplication
$$(a\# b)(c\# d)= ac_R\# b_Rd.$$
$(A,B,R)$ is called a weak smash product structure if $A\#_RB$ is an
associative $k$-algebra with preunit $1_A\# 1_B$. 
The multiplication is associative if and only if
$$R(bd\ot a)=a_{Rr}\ot b_rd_R~~{\rm and}~~
R(b\ot ac)=a_Rc_r\ot b_{Rr},$$
for all $a,c\in A$ and $b,d\in B$. $1_A\# 1_B$ is a preunit if and only
if
$$R(1_B\ot a)=a(1_A)_R\ot (1_B)_R~~{\rm and}~~R(b\ot 1_A)=(1_A)_R\ot (1_B)_Rb.$$

A left-right weak entwining structure is a triple $(A,C,\psi)$, where
$A$ is an algebra, $C$ is a coalgebra, and $\psi:\ A\ot C\to A\ot C$
is a $k$-linear map satisfying the conditions
$$a_\psi\ot \Delta(c^\psi)=a_{\psi\Psi}\ot c_{(1)}^\Psi\ot c_{(2)}^\psi~~;~~
(ab)_\psi\ot c^\psi=a_\psi b_\Psi\ot c^{\Psi\psi};$$
$$1_\psi\ot c^\psi=\varepsilon(c_{(1)}^\psi)1_\psi\ot c_{(2)}~~;~~
a_\psi\varepsilon(c^\psi)=\varepsilon(c^\psi)a1_\psi.$$
Here we use the notation (with summation implicitely understood):
$$\psi(a\ot c)=a_\psi\ot c^\psi.$$
An entwined module is a $k$-module $M$ with a left $A$-action and a right
$C$-coaction such that
$$\rho(am)=a_\psi m_{[0]}\ot m_{[1]}^\psi.$$
The category of entwined modules and left $A$-linear right $C$-colinear
maps is denoted by ${}_A\Mm(\psi)^C$.\\
Let $H$ be a weak bialgebra, and $A$ a right $H$-comodule, which is also an
algebra with unit. $A$ is called a right $H$-comodule algebra if
$\rho(a)\rho(b)=\rho(ab)$ and
$1_{[0]}\ot \varepsilon_t(1_{[1]})=\rho(1)$.\\

From \cite{Bohm}, we recall the following definitions.
Let $C$ be a left $H$-module which is also a coalgebra with counit. $C$
is called a left $H$-comodule algebra if
$\Delta_C(hc)= \Delta_H(h)\Delta_C(c)$ and
\begin{equation}\eqlabel{3.1.1}
\varepsilon_C(hkc)=\varepsilon_H(hk_{(2)})\varepsilon_C(k_{(1)}c),
\end{equation}
for all $c\in C$ and $h,k\in H$. 
Several equivalent definitions are given in 
\cite[Sec. 4]{CDG}. We then call $(H,A,C)$ a 
left-right weak Doi-Hopf datum.
A weak Doi-Hopf module over $(H,A,C)$ is a $k$-module $M$ with a left $A$-action
and a right $C$-coaction, satisfying the following compatibility relation, for
all $m\in M$ and $a\in A$:
\begin{equation}\eqlabel{3.1.2}
\rho(am)=a_{[0]}m_{[0]}\ot a_{[1]}m_{[1]}.
\end{equation}
The category of weak Doi-Hopf modules over $(H,A,C)$ and left $A$-linear
right $C$-colinear maps is denoted by ${}_A\Mm(H)^C$.\\
Let $(H,A,C)$ be a weak left-right Doi-Hopf datum, and consider the map
$$\psi:\ A\ot C\to A\ot C,~~\psi(a\ot c)=a_{[0]}\ot a_{[1]}c.$$
Then $(A,C,\psi)$ is a weak left-right entwining structure, and we have
an isomorphism of categories ${}_A\Mm(H)^C\cong {}_A\Mm(\psi)^C$.\\

Let $(A,C,\psi)$ be a weak left-right entwining structure, and assume that
$C$ is finitely generated projective as a $k$-module, with finite dual basis
$\{(c_i,c_i^*)~|~i=1,\cdots,n\}$. Then we have a weak smash product structure
$(A,C^*,R)$, with $R:\ C^*\ot A\to A\ot C^*$ given by
\begin{equation}\eqlabel{1.41}
R(c^*\ot a)=\hbox{$\sum_i$} \lan c^*,c_i^\psi\ran a_\psi\ot c_i^*.
\end{equation}
We have an isomorphism of categories
\begin{equation}\eqlabel{1.42}
F:\ {}_A\Mm(\psi)^C\to {}_{\ol{A\#_R C^*}}\Mm,
\end{equation}
defined also follows: $F(M)=M$ as a $k$-module, with action
$[a\# c^*]\cdot m=\lan c^*,m_{[1]}\ran am_{[0]}.$ Details can be found
in \cite[Theorem 3.4]{CDG}.

\section{Yetter-Drinfeld modules over weak Hopf algebras}\selabel{2}
Let $H$ be a weak bialgebra. A left-left
Yetter-Drinfeld module is a $k$-module with a left $H$-action and
a left $H$-coaction such that the following conditions hold, for all
$m\in M$ and $h\in H$:
\begin{eqnarray}
&&\hspace*{-2cm}\lambda(m)= m_{[-1]}\ot m_{[0]}\in H\ot_t M;\eqlabel{2.1.1}\\
&&\hspace*{-2cm} h_{(1)}m_{[-1]}\ot h_{(2)}m_{[0]}=(h_{(1)}m)_{[-1]}h_{(2)}\ot 
(h_{(1)}m)_{[0]}.\eqlabel{2.1.2}
\end{eqnarray}
We will now state some equivalent definitions. First we will rewrite the
counit property for Yetter-Drinfeld modules.

\begin{lemma}\lelabel{2.1}
Let $H$ be a weak bialgebra, and 
$\lambda:\ M\to H\ot_t M$, $\rho(m)=m_{[-1]}\ot m_{[0]}$ a $k$-linear map. Then
\begin{equation}\eqlabel{2.1.3}
\varepsilon(m_{[-1]})m_{[0]}=\varepsilon_t(m_{[-1]})m_{[0]}.
\end{equation}
Consequently, in the definition of a Yetter-Drinfeld module, the counit
property $\varepsilon(m_{[-1]})m_{[0]}=m$ can be replaced by
$\varepsilon_t(m_{[-1]})m_{[0]}=m$.
\end{lemma}

\begin{proof}
$$\varepsilon_t(m_{[-1]})m_{[0]}=\varepsilon(1_{(1)}m_{[-1]})1_{(2)}m_{[0]}=
\varepsilon(m_{[-1]})m_{[0]}.$$
\end{proof}

In the case of a weak Hopf algebra,
the compatibility relation \equref{2.1.2} can also be restated:

\begin{proposition}\prlabel{2.2} (cf. \cite[Remark 2.6]{Nenciu})
Let $H$ be a weak Hopf algebra, and $M$ a $k$-module, with a left
$H$-action and a left $H$-coaction. $M$ is a Yetter-Drinfeld module if and only if
\begin{equation}\eqlabel{2.2.1}
\lambda(hm)=h_{(1)}m_{[-1]}S(h_{(3)})\ot h_{(2)}m_{[0]}.
\end{equation}
\end{proposition}

\begin{proof}
Let $M$ be a Yetter-Drinfeld module.
Then we compute
\begin{eqnarray*}
&&\hspace*{-15mm}
h_{(1)}m_{[-1]}S(h_{(3)})\ot h_{(2)}m_{[0]}=
(h_{(1)}m)_{[-1]}h_{(2)}S(h_{(3)})\ot  (h_{(1)}m)_{[0]}\\
&=& (h_{(1)}m)_{[-1]}\varepsilon_t(h_{(2)})\ot  (h_{(1)}m)_{[0]}
\equal{\equref{1.1.3}}
(1_{(1)}hm)_{[-1]}1_{(2)}\ot  (1_{(1)}hm)_{[0]}\\
&\equal{\equref{2.1.2}}& 1_{(1)}(hm)_{[-1]}\ot
1_{(2)}(hm)_{[0]}\equal{\equref{2.1.1}}
(hm)_{[-1]}\ot (hm)_{[0]}=\lambda(hm).
\end{eqnarray*}
Conversely, assume that \equref{2.2.1} holds for all $h\in H$ and $m\in M$.
Taking $h=1$ in \equref{2.2.1}, we find
\begin{eqnarray*}
\lambda(m)&=&1_{(1)}m_{[-1]}S(1_{(3)})\ot 1_{(2)}m_{[0]}\\
&=& 1_{(1)}m_{[-1]}S(1_{(2')})\ot 1_{(2)}1_{(1')}m_{[0]}\in H\ot_t M
\end{eqnarray*}
and
\begin{eqnarray}
\lambda(m)&=& 1_{(1)}m_{[-1]}S(1_{(3)})\ot 1_{(2)}m_{[0]}
= 1_{(1)}m_{[-1]}S(1_{(2')})\ot 1_{(1')}1_{(2)}m_{[0]}\nonumber\\
&=& m_{[-1]}S(1_{(2')})\ot 1_{(1')}m_{[0]}.\eqlabel{2.2.2}
\end{eqnarray}
Now
\begin{eqnarray*}
&&\hspace*{-15mm}
(h_{(1)}m)_{[-1]}h_{(2)}\ot (h_{(1)}m)_{[0]}\equal{\equref{2.2.1}}
h_{(1)}m_{[-1]}S(h_{(3)})h_{(4)}\ot h_{(2)}m_{[0]}\\
&=& h_{(1)}m_{[-1]}\varepsilon_s(h_{(3)})\ot h_{(2)}m_{[0]}
\equal{\equref{1.6.3}}h_{(1)}m_{[-1]}S(1_{(2)})\ot h_{(2)}1_{(1)}m_{[0]}\\
&\equal{\equref{2.2.2}}&
h_{(1)}m_{[-1]}\ot h_{(2)}m_{[0]},
\end{eqnarray*}
as needed.
\end{proof}

\begin{corollary}\colabel{2.2bis}
Let $M$ be a left-left Yetter-Drinfeld module. For all $y\in H_s$, $z\in H_t$
and $m\in M$, we have
\begin{equation}\eqlabel{2.2bis.1}
\lambda(zm)=zm_{[-1]}\ot m_{[0]}~~;~~\lambda(ym)=m_{[-1]}S(y)\ot m_{[0]}.
\end{equation}
\end{corollary}

\begin{proof}
\begin{eqnarray*}
\lambda(zm)&\equal{(\ref{eq:1.1.4},\ref{eq:2.2.1})}&
1_{(1)}zm_{[-1]}S(1_{(3)})\ot 1_{(2)}m_{[0]}\\
&\equal{\equref{1.1.6}}&
z1_{(1)}m_{[-1]}S(1_{(3)})\ot 1_{(2)}m_{[0]}
\equal{\equref{2.2.1}} zm_{[-1]}\ot m_{[0]}.
\end{eqnarray*}
The other assertion is proved in a similar way.
\end{proof}

\begin{corollary}\colabel{2.3}
Let $M$ be a left-left Yetter-Drinfeld module over a weak Hopf algebra
with bijective antipode.
Then we have the following identities, for all $m\in M$:
\begin{equation}\eqlabel{2.3.1}
1_{(1)}m_{[0]}\ot 1_{(2)}S^{-1}(m_{[-1]})= m_{[0]}\ot S^{-1}(m_{[-1]});
\end{equation}
\begin{equation}\eqlabel{2.3.2}
\varepsilon_s(S^{-2}(m_{[-1]}))m_{[0]}=m.
\end{equation}
\end{corollary}

\begin{proof}
Apply $S^{-1}$ to the first factor of \equref{2.2.2}, and then switch the
two tensor factors. Then we obtain \equref{2.3.1}. \equref{2.3.2} is proved
as follows:
\begin{eqnarray*}
&&\hspace*{-2cm}
m=\varepsilon_t(m_{[-1]})m_{[0]}\equal{\equref{2.1.3}}
\varepsilon(m_{[-1]})m_{[0]}
= \varepsilon(S^{-1}(m_{[-1]}))m_{[0]}\\
&\equal{\equref{2.3.1}}&
\varepsilon(1_{(2)}S^{-1}(m_{[-1]}))1_{(1)}m_{[0]}
\equal{\equref{1.6.2}}\varepsilon_s(S^{-2}(m_{[-1]}))m_{[0]}.
\end{eqnarray*}
\end{proof}

The category of left-left Yetter-Drinfeld modules and left $H$-linear,
left $H$-colinear maps will be denoted by ${}_H^H\YD$.

\begin{example}\exlabel{2.3}
Let $G$ be a groupoid, and $kG$ the corresponding groupoid algebra.
Then $kG$ is a weak Hopf algebra. Let $M$ be a left-left Yetter-Drinfeld
module. Then $M$ is a $kG$-comodule, so $M$ is graded by the set $G$,
that is
$$M=\bigoplus_{\sigma\in G_1} M_\sigma,$$
and $\lambda(m)=\sigma\ot m$ if and only if $m\in M_\sigma$, or
${\rm deg}(m)=\sigma$.\\
Recall that the unit element of $kG$ is $1=\sum_{x\in G_0} x$, where
$x$ is the identity morphism of the object $x\in G_0$. 
Take $m\in M_\sigma$. Using \equref{2.2.1},
we find
$$\lambda(m)=\lambda(1m)=\sum_{x\in G} x\sigma x\ot xm=0,$$
unless $s(\sigma)=\tau(\sigma)=x$. So we have
$$M= \bigoplus_{\sigma\in G_1\atop s(\sigma)=t(\sigma)} M_\sigma.$$
Take $m\in M_\sigma$, with $s(\sigma)=\tau(\sigma)$, and $\tau\in G_1$.
It follows from \equref{2.2.1} that
$\lambda(\tau m)= \tau\sigma\tau^{-1}\ot \tau m =0$,
unless $s(\tau)=x$. If $s(\tau)=x$, then ${\rm deg}(\tau m)=\tau\sigma\tau^{-1}$.
\end{example}

\begin{theorem}\thlabel{2.4}
Let $H$ be a weak bialgebra. Then the
category ${}_H^H\YD$ is isomorphic to the weak left center $\Ww_l({}_H\Mm)$
of the category of left $H$-modules. If $H$ is a weak Hopf algebra with bijective antipode,
then ${}_H^H\YD$ is isomorphic to the left center $\Zz_l({}_H\Mm)$
\end{theorem}

\begin{proof}
We will restrict to a brief description of the connecting functors; for more detail
(in the left-right case),
we refer to \cite[Lemma 4.3]{Nenciu}.
Take $(M,\sigma_{M,-})\in \Ww_l({}_H\Mm)$. For each left $H$-module $V$,
we have a map $\sigma_{M,V}:\ M\ot_t V\to V\ot_t M$ in ${}_H\Mm$.
We will show that the map
$$\lambda:\ M\to H\ot_t M,~~\lambda(m)=\sigma_{M,H}(1_{(1)}m\ot 1_{(2)})=
m_{[-1]}\ot m_{[0]}$$
makes $M$ into a Yetter-Drinfeld module. Conversely, let $(M,\lambda)$ is a Yetter-Drinfeld module; a natural
transformation $\sigma$ is then defined by the formula
\begin{equation}\eqlabel{2.4.1}
\sigma_{M,V}(1_{(1)}m\ot 1_{(2)}v)=m_{[-1]}v\ot m_{[0]}.
\end{equation}
Straightforward
computations show that $(M,\sigma)\in \Ww_l({}_H\Mm)$. If $H$
is a Hopf algebra with invertible antipode, then the inverse of
$\sigma_{M,V}$ is 
\begin{equation}\eqlabel{2.4.2}
\sigma_{M,V}^{-1}(1_{(1)}v\ot 1_{(2)}m)=m_{[0]}\ot S^{-1}(m_{[-1]})v.
\end{equation}
\end{proof}

From now on, we assume that $H$ is a weak Hopf algebra with bijective antipode.
Since the left center of a monoidal category is a braided monoidal category,
it follows from \thref{2.4} that ${}_H^H\YD$ is a braided monoidal category;
a direct but long proof can be given: see \cite[Prop. 2.7]{Nenciu}.
The monoidal structure can be computed using \equref{center2}.
Take $M,N\in {}_H^H\YD$, the $H$-coaction on $M\ot_tN$ is given by the 
formula
$$\lambda(1_{(1)}m\ot 1_{(2)}n)=
((\sigma_{M,H}\ot N)\circ (M\ot \sigma_{N,H}))
(1_{(1')}(1_{(1)}m\ot 1_{(2)}n)\ot 1_{(2')}).$$
Observe that
\begin{eqnarray*}
x&=& 1_{(1')}(1_{(1)}m\ot 1_{(2)}n)\ot 1_{(2')}\\
&=& 1_{(1')}1_{(1)}m\ot 1_{(1'')}1_{(2')}1_{(2)}n\ot 1_{(2'')}=
1_{(1)}m\ot 1_{(1'')}1_{(2)}n\ot 1_{(2'')},
\end{eqnarray*}
so that
\begin{eqnarray*}
&&\hspace*{-2cm}
(M\ot \sigma_{N,H})(x)=
1_{(1)}m\ot (1_{(2)}n)_{[-1]}\ot (1_{(2)}n)_{[0]}\\
&=& 1_{(1)}m\ot 1_{(2)}n_{[-1]}S(1_{(4)})\ot 1_{(3)}n_{[0]}\\
&=& 1_{(1)}m\ot 1_{(2)}1_{(1')}n_{[-1]}S(1_{(3')})\ot 1_{(2')}n_{[0]}\\
&=& 1_{(1)}m\ot 1_{(2)}n_{[-1]}\ot n_{[0]}
\end{eqnarray*}
and
\begin{equation}\eqlabel{2.5.1}
\lambda(1_{(1)}m\ot 1_{(2)}n)=m_{[-1]}n_{[-1]}\ot m_{[0]}\ot n_{[0]}.
\end{equation}
We compute the left $H$-coaction on $H_t$ using \equref{center2bis} and 
\equref{2.4.1}. For any $z\in H_t$, this gives
\begin{eqnarray}
\lambda(z)&=& \sigma_{H_t,H}((1_{(1)}\leftact z)\ot 1_{(2)})=
 r_M^{-1}(l_M((1_{(1)}\leftact z)\ot 1_{(2)}))\nonumber\\
&=&
r_M^{-1}(z)= 1_{(1)}z \ot 1_{(2)}=\Delta(z).\eqlabel{2.5.2}
\end{eqnarray}
The braiding and its inverse are given by the formulas
$$\sigma_{M,N}(1_{(1)} m\ot 1_{(2)}n)=m_{[-1]}n\ot m_{[0]}~~;~~
\sigma^{-1}_{M,N}(1_{(1)} n\ot 1_{(2)}m)=
m_{[0]}\ot S^{-1}(m_{[-1]})n.$$

A left-right
Yetter-Drinfeld module is a $k$-module with a left $H$-action and
a right $H$-coaction such that the following conditions hold, for all
$m\in M$ and $h\in H$:
\begin{eqnarray}
&&\hspace*{-2cm}\rho(m)= m_{[0]}\ot m_{[1]}\in M\ot_t H;\eqlabel{2.6.1}\\
&&\hspace*{-2cm} h_{(1)}m_{[0]}\ot h_{(2)}m_{[1]}=(h_{(2)}m)_{[0]}\ot 
(h_{(2)}m)_{[1]}h_{(1)}.\eqlabel{2.6.2}
\end{eqnarray}
The category of left-right Yetter-Drinfeld modules and left $H$-linear
right $H$-colinear maps is denoted by ${}_H\YD^H$.

\begin{proposition}\prlabel{2.6}
Let $H$ be a weak Hopf algebra with bijective antipode. Then the
category ${}_H\YD^H$ is isomorphic to the right center $\Zz_r({}_H\Mm)$.
\end{proposition}

\begin{proof}
Take $(M, \tau_{-,M})\in \Zz_r({}_H\Mm)$. We know from \prref{1.0.1}
that $(M, \sigma_{M,-}=\tau^{-1}_{-,M})\in \Zz_l({}_H\Mm)$. Take the
corresponding left-left Yetter-Drinfeld $(M,\lambda)$, as in
\thref{2.4}, and define $\rho:\ M\to M\ot H$ by
\begin{equation}\eqlabel{2.6.0}
\rho(m)=m_{[0]}\ot m_{[1]}=m_{[0]}\ot S^{-1}(m_{[-1]}).
\end{equation}
It follows from \equref{2.3.1} that $\rho(m)\in M\ot_t H$. 
The coassociativity of $\rho$ follows immediately from the coassociativity
of $\lambda$ and the anti-comultiplicativity of $S^{-1}$. Also
$$\varepsilon(m_{[1]})m_{[0]}=\varepsilon(S^{-1}(m_{[-1]}))m_{[0]}
=\varepsilon(m_{[-1]})m_{[0]}=m.$$
From
\equref{2.4.2}, it follows that
\begin{equation}\eqlabel{2.6.3}
\tau_{V,M}(1_{(1)}v\ot 1_{(2)}m)=m_{[0]}\ot m_{[1]}v.
\end{equation}
In particular, $\tau_{M,H}(1_{(1)}\ot 1_{(2)}m)=\rho(m)$, and the
fact that $\tau_{M,H}$ is left $H$-linear implies \equref{2.6.2}.
Hence $(M,\rho)$ is a left-right Yetter-Drinfeld module.\\
Conversely, if  $(M,\rho)$ is a left-right Yetter-Drinfeld module, then
$(M,\tau_{-,M})$, with $\tau$ defined by \equref{2.6.3} is an object
of $\Zz_r({}_H\Mm)$.
\end{proof}

\begin{corollary}\colabel{2.7}
Let $M$ be a $k$-module with a left $H$-action and a right $H$-coaction.
Then $M$ is a 
left-right Yetter-Drinfeld module if and only if
\begin{equation}\eqlabel{2.7.1}
\rho(hm)=h_{(2)}m_{[0]}\ot h_{(3)}m_{[1]}S^{-1}(h_{(1)}).
\end{equation}
\end{corollary}

\begin{corollary}\colabel{2.7a}
Let $M$ be a left-right Yetter-Drinfeld module. For all $y\in H_s$,
$z\in H_t$ and $m\in M$, we have that
\begin{equation}\eqlabel{2.7a.1}
\rho(ym)=m_{[0]}\ot ym_{[1]}~~;~~\rho(zm)=m_{[0]}\ot m_{[1]}S^{-1}(z).
\end{equation}
\end{corollary}

\begin{corollary}\colabel{2.7b}
Let $M$ be a left-right Yetter-Drinfeld module. Then 
\begin{equation}\eqlabel{2.7b.1}
1_{(2)}m_{[0]}\ot m_{[1]}S^{-1}(1_{(1)})=\rho(m),
\end{equation}
for all $m\in M$.
\end{corollary}

\begin{proof}
Apply $S^{-1}\ot M$ to 
$\lambda(m)=1_{(1)}S(m_{[1]})\ot 1_{(2)}m_{[0]}$.
\end{proof}

\begin{corollary}\colabel{2.8}
The category ${}_H\YD^H$ is a braided monoidal category, isomorphic to
${}_H^H\YD^{\rm in}$.
\end{corollary}

In a similar way, we can introduce right-right and right-left
Yetter-Drinfeld modules. The categories $\YD_H^H$ and ${}^H\YD_H$
of right-right and right-left
Yetter-Drinfeld modules are isomorphic to the right and left center
of $\Mm_H$. Let us summarize the results.\\
A right-right Yetter-Drinfeld module is a $k$-module $M$ with a
right $H$-action and a right $H$-coaction such that
\begin{eqnarray}
&&\hspace*{-2cm}\rho(m)= m_{[0]}\ot m_{[1]}\in M\ot_s H;\eqlabel{2.9.1}\\
&&\hspace*{-2cm} m_{[0]}h_{(1)}\ot m_{[1]}h_{(2)}=(mh_{(2)})_{[0]}\ot 
h_{(1)}(mh_{(2)})_{[1]};\eqlabel{2.9.2}
\end{eqnarray}
or, equivalently,
\begin{equation}\eqlabel{2.9.3}
\rho(mh)=m_{[0]}h_{(2)}\ot S(h_{(1)})m_{[1]}h_{(3)}.
\end{equation}
The counit condition $m=\varepsilon(m_{[1]})m_{[0]}$ is equivalent to
$$m= m_{[0]}\varepsilon(m_{[1]}).$$
The natural isomorphism $\tau_{-,M}$ corresponding to $(M,\rho)\in \YD_H^H$
and its inverse are given by the formulas
\begin{equation}\eqlabel{2.9.4}
\tau_{M,V}(v1_{(1)}\ot m1_{(2)})=m_{[0]}\ot mv_{[1]}~~;~~
\tau^{-1}_{M,V}(m1_{(1)}\ot v1_{(2)})=vS^{-1}(m_{[1]})\ot m_{[0]}.
\end{equation}
Furthermore
$$m_{[0]}\varepsilon_t(S^{-2}(m_{[1]}))=m,$$
and $S^{-1}(m_{[1]})\ot m_{[0]}\in H\ot_s M$.\\
The monoidal structure on $\YD_H^H$ is given by the formula
$$\rho(m1_{(1)}\ot n1_{(2)})=m_{[0]}\ot n_{[0]}\ot m_{[1]}n_{[1]}.$$
The braiding is given by \equref{2.9.4}. The category $\YD_H^H$
is isomorphic as a braided monoidal category to $\Zz_r(\Mm_H)$.\\

Let $M$ be a right $H$-module and a left $H$-comodule. $M$ is a right-left
Yetter-Drinfeld module if one of the three following equivalent conditions
is satisfied, for all $m\in M$ and $h\in H$:\\
1) $\lambda(m)\in H\ot_sM$ and
$$h_{(2)}(mh_{(1)})_{[0]}\ot (mh_{(1)})_{[1]}=
m_{[-1]}h_{(1)}\ot m_{[0]}h_{(2)},$$
2)
$\lambda(mh)=S^{-1}(h_{(3)})m_{[-1]}h_{(1)}\ot m_{[0]}h_{(2)}$;\\
3) $(M,\rho)$, with $\rho(m)=m_{[0]}\ot S(m_{[-1]})$ is a right-right
Yetter-Drinfeld module.\\

The category of right-left Yetter-Drinfeld modules, ${}^H\YD_H$, is
a braided monoidal category. The monoidal structure and the braiding are given by
$$\lambda(m1_{(1)}\ot n1_{(2)})=m_{[-1]}n_{[-1]}\ot m_{[0]}\ot n_{[0]};$$
$$\sigma_{M,N}(m1_{(1)}\ot n1_{(2)})=nm_{[-1]}\ot m_{[0]}.$$
As a braided monoidal category, ${}^H\YD_H$ is isomorphic to
$\Zz_l(\Mm_H)$ and $(\YD_H^H)^{\rm in}$.\\

The antipode $S:\ H\to H^{\rm op,cop}$ is an isomorphism of weak Hopf algebras.
Observe that the target map of $H^{\rm op,cop}$ is $\varepsilon_s$,
and that its source map is $\varepsilon_t$. Thus $S$ induces an isomorphism
between the monoidal categories ${}_H\Mm$ and ${}_{H^{\rm op,cop}}\Mm$.
We also have a monoidal isomorphism $F:\ {}_{H^{\rm op,cop}}\Mm\to \ol{\Mm}_H$,
given by
$$F(M)=M,~~mh=h^{\rm op,cop}m.$$
indeed, in ${}_{H^{\rm op,cop}}\Mm$, $M\ot_t N$ is generated by elements of
the form $1_{(2)}m\ot 1_{(1)}n$, and $F(M\ot_t N)$ is generated by elements of
the form $m1_{(2)}\ot n1_{(1)}$. $F(N)\ot_sF(M)$ is generated by elements of
the form $n1_{(1)}\ot m1_{(2)}$, and it follows that the switch map is an
isomorphism $F(M\ot_t N)\to F(N)\ot_s F(M)$. We conclude from
\prref{1.0.2} that we have isomorphisms of braided monoidal categories
$${}_H^H\YD\cong \Zz_l({}_H\Mm)\cong \Zz_l({}_{H^{\rm op,cop}}\Mm)
\cong \Zz_l(\ol{\Mm}_H)\cong \ol{\Zz_r(\Mm_H)}\cong \ol{\YD}_H^H.$$
This isomorphism can be described explicitely as follows:
$$F:\ {}_H^H\YD\to \ol{\YD}_H^H,~~F(M)=M,$$
with
$$m\cdot h=S^{-1}(h)m~~;~~\rho(m)=m_{[0]}\ot S(m_{[-1]}).$$
We summarize our results as follows:

\begin{theorem}\thlabel{2.10}
Let $H$ be a weak Hopf algebra with bijective antipode. Then we have the following isomorphisms of braided monoidal categories:
$$
{}_H^H\YD\cong {{}_H\YD^H}^{\rm in}\cong \ol{\YD _H^H}
\cong \ol{{}^H\YD_H}^{\rm in}.$$
\end{theorem}

\section{Yetter-Drinfeld modules are Doi-Hopf modules}\selabel{3}
It was shown in \cite{cmz} that Yetter-Drinfeld modules (over a classical Hopf
algebra) can be considered as Doi-Hopf modules, and, a fortiori, as
entwined modules, and as comodules over a coring (see \cite{BW}).
Weak Doi-Hopf modules were introduced by B\"ohm
\cite{Bohm}, and they are special cases of weak entwined modules
(see \cite{CDG}), and these are in turn examples of comodules over a
coring (see \cite{BW}). In this Section, we will show that Yetter-Drinfeld modules
over weak Hopf algebras are special cases of weak Doi-Hopf modules. We will discuss
the left-right case.

\begin{proposition}\prlabel{3.1}
Let $H$ be a weak Hopf algebra with a bijective antipode. Then $H$ is a
right $H\ot H^{\rm op}$-comodule algebra, with $H$-coaction
$$\rho(h)=h_{(2)}\ot S^{-1}(h_{(1)})\ot h_{(3)}.$$
\end{proposition}

\begin{proof}
It is easy to verify that $H$ is  a right $H\ot H^{\rm op}$-comodule and that
$\rho(hk)=\rho(h)\rho(k)$. 
Recall that $H_t=\im(\varepsilon_t)=\im(\ol{\varepsilon}_t)$.
The target map of $H^{\rm op}\ot H$ is $\ol{\varepsilon}_t\ot\varepsilon_t$.
We now have
\begin{eqnarray*}
&&\hspace*{-2cm}
1_{[0]}\ot (\ol{\varepsilon}_t\ot\varepsilon_t)(1_{[1]})=
1_{(2)}1_{(1')}\ot \ol{\varepsilon}_t(S^{-1}(1_{(1)}))\ot 
\varepsilon_t(1_{(2')})\\
&=& 1_{(2)}1_{(1')}\ot S^{-1}(1_{(1)})\ot 1_{(2')}=\rho(1),
\end{eqnarray*}
where we used the fact that $S^{-1}(1_{(1)})\ot 1_{(2)}\in H_t\ot H_t$.
\end{proof}

\begin{proposition}\prlabel{3.2}
Let $H$ be a weak Hopf algebra with a bijective antipode. Then $H$ is a
left $H^{\rm op}\ot H$-module coalgebra with left action
$$(k\ot h)\triangleright c= hck.$$
\end{proposition}

\begin{proof}
We easily compute that
\begin{eqnarray*}
&&\hspace*{-2cm}
\varepsilon((m\ot l)(k_{(2)}\ot h_{(2)}))\varepsilon((k_{(1)}\ot h_{(1)})\triangleright c)\\
&=& \varepsilon(k_{(2)}m)\varepsilon(lh_{(2)})\varepsilon(h_{(1)}ck_{(1)})\\
&=& \varepsilon(lhckm)=\varepsilon(((m\ot l)(k\ot h))\triangleright c).
\end{eqnarray*}
The other conditions are easily verified.
\end{proof}

\begin{corollary}\colabel{3.3}
Let $H$ be a weak Hopf algebra with bijective antipode. Then we have
a weak Doi-Hopf datum $(H^{\rm op}\ot H, H,H)$ and the categories
${}_H\Mm(H^{\rm op}\ot H)^H$ and ${}_H\YD^H$ are isomorphic.
\end{corollary}

\begin{proof}
The compatibility relation \equref{3.1.2} reduces to \equref{2.7.1}.
\end{proof}

As we have seen in \seref{1.4}, weak Doi-Hopf modules are special cases
of entwined modules. The entwining map $\psi:\ H\ot H\to H\ot H$
corresponding to the weak Doi-Hopf datum $(H^{\rm op}\ot H, H,H)$ is given
by
\begin{equation}\eqlabel{3.3.1}
\psi(h\ot k)=h_{(2)}\ot h_{(3)}kS^{-1}(h_{(1)}).
\end{equation}

\section{The Drinfeld double}\selabel{4}
Now we consider the particular case where $H$ is finitely generated and
projective as a $k$-module, with finite dual basis
$\{(h_i,h_i^*)~|~i=1,\cdots,n\}$. Then $H^*$ is also a weak Hopf algebra,
in view of the selfduality of the axioms of a weak Hopf algebra.
Recall that the comultiplication is given by the formula
$\lan \Delta(h^*),h\ot k\ran =\lan h^*,hk\ran$; the counit is evaluation at $1$.
Also recall that $H^*$ is an $H$-bimodule, with left and right $H$-action
$$\lan h\leftact h^*\rightact k, l\ran =\lan h^*, klh\ran,$$
or
\begin{equation}\eqlabel{4.0}
h\leftact h^*\rightact k=\lan h_{(1)}^*,k\ran \lan h_{(3)}^*,h\ran h_{(2)}^*.
\end{equation}
Using \equref{1.41}, we find a weak smash
product structure $(H,H^*,R)$, with $R:\ H^*\ot H\to H\ot H^*$ given by
\begin{eqnarray}
R(h^*\ot h)&=&
\hbox{$\sum_i$} \lan h^*, h_{(3)}h_iS^{-1}(h_{(1)})\ran h_{(2)}\ot h_i^*\nonumber\\
&=& \hbox{$\sum_i$} \lan S^{-1}(h_{(1)})\leftact h^*\rightact h_{(3)}, h_i\ran h_{(2)}\ot h_i^*\nonumber\\
&=& h_{(2)} \ot \Bigl(S^{-1}(h_{(1)})\leftact h^*\rightact h_{(3)}\Bigr).
\eqlabel{4.1}
\end{eqnarray}

From \seref{1.4}, we know that $H\#_R H^*$, which we will also denote
by $H\bowtie H^*$, is an associative algebra with preunit $1\# \varepsilon$.
Using \equref{smash}, we compute the multiplication rule on $H\bowtie H^*$.
\begin{eqnarray}
&&\hspace*{-2cm} (h\bowtie h^*)(k\bowtie k^*)=
\hbox{$\sum_i$} \lan h^*, k_{(3)}h_iS^{-1}(k_{(1)})\ran hk_{(2)}
\bowtie h_i^**k^*\nonumber\\
&=& hk_{(2)}\bowtie (S^{-1}(k_{(1)})\leftact h^*\rightact k_{(3)})*k^*
\eqlabel{4.2}\\
&=& hk_{(2)} \bowtie \lan h_{(1)}^*, k_{(3)}\ran \lan h_{(3)}^*,S^{-1}(k_{(1)})\ran
h_{(2)}^**k^*\eqlabel{4.3}.
\end{eqnarray}
We have a projection $p:\ H\bowtie H^*\to H\bowtie H^*$,
$$p(h\bowtie h^*)=(1\bowtie \varepsilon)(h\bowtie h^*)=(h\bowtie h^*)(1\bowtie \varepsilon)= (h\bowtie h^*)(1\bowtie \varepsilon)^2,$$
and $D(H)=\ol{H\bowtie H^*}= (H\bowtie H^*)/\Ker p$ is a $k$-algebra with
unit $[1\bowtie \varepsilon]$, which we call the {\sl Drinfeld double}
of $H$. $D(H)$ is also isomorphic to $\ul{H\bowtie H^*}=\im(p)$, which is
a $k$-algebra with unit $(1\bowtie \varepsilon)^2$. Observe that the multiplication rule
\equref{4.3} is the same as in \cite{Bohm,Nenciu}. We show that the ideal
$J$ that is divided out in \cite{Bohm,Nenciu} is equal to $\Ker p$, and this will
imply that $D(H)$ is equal to the Drinfeld double introduced in \cite{Bohm,Nenciu}.
We first need some Lemmas.

\begin{lemma}\lelabel{4.1}
Let $H$ a weak bialgebra. For all $h^*\in H^*$, $y\in H_s$ and $z\in H_t$,
we have
\begin{eqnarray}
h^**(y\leftact\varepsilon)&=&
\lan h_{(2)}^*, y\ran  h_{(1)}^*=y\leftact h^*\eqlabel{4.1.1}\\
h^**(\varepsilon\rightact y)&=&
\lan h_{(1)}^*, y\ran  h_{(2)}^*=h^*\rightact y\eqlabel{4.1.2}\\
(z\leftact \varepsilon)*h^*&=&
\lan h_{(2)}^*,z\ran h_{(1)}^*=z\leftact h^*\eqlabel{4.1.3}\\
(\varepsilon\rightact z)*h^*&=&
\lan h_{(1)}^*,z\ran h_{(2)}^*=h^*\rightact z\eqlabel{4.1.4}
\end{eqnarray}
\end{lemma}

\begin{proof}
We only prove \equref{4.1.3}. For all $h\in H$, we have
\begin{eqnarray*}
&&\hspace*{-2cm}
\lan (z\leftact \varepsilon)*h^*, h\ran=
\lan \varepsilon, h_{(1)}z\ran \lan h^*,h_{(2)}\ran
= \lan \varepsilon, h_{(1)}1_{(1)}z\ran \lan h^*,h_{(2)}1_{(2)}\ran\\
&=& \lan \varepsilon*h^*, hz\ran
= \lan h^*, hz\ran=\lan z\leftact h^*,h\ran = \lan h_{(2)}^*,z\ran
\lan h_{(1)}^*,h\ran.
\end{eqnarray*}
\end{proof}

\begin{lemma}\lelabel{4.1a}
Let $H$ be a weak Hopf algebra with bijective antipode. For all $y\in H_s$,
$z\in H_t$, we have
\begin{equation}\eqlabel{4.1.5}
S^{-1}(z)\leftact \varepsilon=z\leftact\varepsilon~~{\rm and}~~
\varepsilon\rightact y=\varepsilon\rightact S^{-1}(y).
\end{equation}
\end{lemma}

\begin{proof}
For all $h\in H$, we have
\begin{eqnarray*}
&&\hspace*{-2cm}
\lan S^{-1}(z)\leftact \varepsilon,h\ran= \varepsilon(hS^{-1}(z))\equal{\equref{1.1.2}}
\varepsilon(h1_{(1)})\varepsilon(1_{(2)}S^{-1}(z))
\equal{\ref{le:1.5}}
\varepsilon(h1_{(1)})\varepsilon(zS(1_{(2)}))\\
&\equal{\equref{1.6.3}}&
\varepsilon(h_{(1)})\varepsilon(z\varepsilon_s(h_{(2)}))
\equal{\equref{1.1.6}}\varepsilon(\varepsilon_s(h)z)
\equal{\equref{1.1.3a}} \varepsilon(hz)=\lan z\leftact\varepsilon, h\ran.
\end{eqnarray*}
The second statement can be proved in a similar way.
\end{proof}

\begin{proposition}\prlabel{4.2}
Let $H$ be a finitely generated projective weak Hopf algebra. Then $\Ker(p)$
is the $k$-linear span $J$ of elements of the form
$$A=hz\bowtie h^*-h\bowtie (z\leftact \varepsilon)*h^*~~{\rm and}~~
B=hy\bowtie h^*-h\bowtie (\varepsilon\rightact y)*h^*,$$
where $h\in H$, $h^*\in H^*$, $y\in H_s$ and $z\in H_t$.
\end{proposition}

\begin{proof}
$A\in \Ker(p)$ since
\begin{eqnarray*}
&&\hspace*{-2cm}
(1\bowtie \varepsilon)(hz\bowtie h^*)\equal{\equref{4.3}}
h_{(2)}1_{(2)}\bowtie \varepsilon_{(2)}*h^*
\lan \varepsilon_{(1)},h_{(3)}1_{(3)}\ran
\lan \varepsilon_{(3)},S^{-1}(h_{(1)}1_{(1)}z)\ran\\
&=& h_{(2)}\bowtie \varepsilon_{(2)}*h^*
\lan \varepsilon_{(1)},h_{(3)}\ran
\lan \varepsilon_{(3)},S^{-1}(z)\ran
\lan \varepsilon_{(4)},S^{-1}(h_{(1)})\ran\\
&\equal{\equref{4.1.1}}& h_{(2)}\bowtie (\varepsilon_{(2)}*
(S^{-1}(z)\leftact \varepsilon)*h^*)
\lan \varepsilon_{(1)},h_{(3)}\ran
\lan \varepsilon_{(3)},S^{-1}(h_{(1)})\ran\\
&\equal{\equref{4.3}}& (1\bowtie \varepsilon)
\bigl(h\bowtie ((S^{-1}(z)\leftact \varepsilon)*h^*)\bigr)
\equal{\equref{4.1.5}} (1\bowtie \varepsilon)(h\bowtie (z\leftact\varepsilon)*h^*).
\end{eqnarray*}
In a similar way, $B\in \Ker(p)$:
\begin{eqnarray*}
&&\hspace*{-2cm}
(1\bowtie \varepsilon)(hy\bowtie h^*)\equal{\equref{4.3}}
(h_{(2)}1_{(2)}\bowtie \varepsilon_{(2)}^**h^*)
\lan \varepsilon_{(1)},h_{(3)}1_{(3)}y\ran
\lan\varepsilon_{(3)}, S^{-1}(h_{(1)}1_{(1)})\ran\\
&=& (h_{(2)}\bowtie \varepsilon_{(3)}^**h^*)
\lan \varepsilon_{(1)},h_{(3)}\ran
\lan \varepsilon_{(2)},y\ran
\lan\varepsilon_{(4)}, S^{-1}(h_{(1)})\ran\\
&\equal{\equref{4.1.2}}&
(h_{(2)}\bowtie (\varepsilon_{(2)}\rightact y)*h^*)
\lan \varepsilon_{(1)},h_{(3)}\ran\lan\varepsilon_{(3)}, S^{-1}(h_{(1)})\ran\\
&\equal{\equref{4.3}}& (1\bowtie \varepsilon)
(h\bowtie (\varepsilon\rightact y)*h^*).
\end{eqnarray*}
This shows that $J\subset \Ker(p)$. We now compute for all $h\in H$ and
$h^*\in H^*$ that
$$
(h\bowtie h^*)(1\bowtie\varepsilon)\equal{\equref{4.3}}
(h1_{(2)}1_{(1')}\bowtie h_{(2)}^*)
\lan h_{(1)}^*,1_{(2')}\ran \lan h_{(3)}^*, S^{-1}(1_{(1)})\ran,$$
and
\begin{eqnarray*}
&&\hspace*{-20mm}
\Bigl(h\bowtie (S^{-1}(1_{(2)})\leftact\varepsilon)*(\varepsilon\rightact 1_{(1')})
*h_{(2)}^*\Bigr) \lan h_{(1)}^*,1_{(2')}\ran \lan h_{(3)}^*,S^{-1}(1_{(1)})\ran\\
&\equal{\equref{4.0}}&
\Bigl(h\bowtie \varepsilon_{(1)}*\varepsilon_{(2')}*h_{(2)}^*\Bigr)
\lan \varepsilon_{(2)}, S^{-1}(1_{(2)})\ran \lan \varepsilon_{(1')}, 1_{(1')}\ran\\
&&\hspace*{1cm}
\lan h_{(1)}^*,1_{(2')}\ran \lan h_{(3)}^*,S^{-1}(1_{(1)})\ran\\
&=& \Bigl(h\bowtie \varepsilon_{(1)}*\varepsilon_{(2')}*h_{(2)}^*\Bigr)
\lan \varepsilon_{(1')}*h_{(1)}^*, 1\ran
\lan \varepsilon_{(2)}*h_{(3)}^*, S^{-1}(1)\ran\\
&=& \Bigl(h\bowtie \varepsilon_{(1)}*h_{(1)}^*\Bigr)
\lan \varepsilon_{(2)}*h_{(2)}^*, 1\ran=
 h\bowtie (\varepsilon* h^*)=h\bowtie h^*.
\end{eqnarray*}
Observing that
\begin{eqnarray*}
&&\hspace*{-2cm}
hzy\bowtie h^* - 
h\bowtie ((S^{-1}(z)\leftact \varepsilon)*(\varepsilon\rightact y)*h^*)\\
&=& hzy\bowtie h^* - hz\bowtie (\varepsilon\rightact y) *h^*)\\
&+&hz\bowtie (\varepsilon\rightact y) *h^*)
- h\bowtie ((S^{-1}(z)\leftact \varepsilon)*(\varepsilon\rightact y)*h^*)\in J,
\end{eqnarray*}
it follows that $(h\bowtie h^*)(1\bowtie\varepsilon)- (h\bowtie h^*)\in J$,
for all $h\in H$ and $h^*\in H^*$. If $x\in \Ker(p)$, then
$x(1\bowtie \varepsilon)=0$, and $x= x-x(1\bowtie \varepsilon)\in J$.
We conclude that $\Ker(p)\subset J$, finishing our proof.
\end{proof}

We now recall the following results from \cite{NVT}. On $H^*\ot H$, there
exists an associative multiplication
\begin{eqnarray*}
(h^*\ot h)(k^*\ot k)&=&
k_{(2)}^*h^*\ot h_{(2)}k \lan S(h_{(1)}),k_{(1)}^*\ran
\lan h_{(3)}, k_{(3)}^*\ran\\
&=& (h_{(3)}\leftact k^*\rightact S(h_{(1)}))*h^*\ot h_{(2)}k.
\end{eqnarray*}
The $k$-module $I$ generated by elements of the form
$$A'=h^*\ot hz-(\varepsilon\rightact z)h^*\ot h~~{\rm and}~~
B'=h^*\ot yh-(y\leftact\varepsilon)h^*\ot h$$
is a two-sided ideal of $H^*\ot H$. The quotient $D'(H)=(H^*\ot H)/I$
is an algebra with unit element $\varepsilon\ot 1$. It is a weak Hopf algebra,
with the following comultiplication, counit and antipode:
\begin{eqnarray}
&&\hspace*{-2cm} \Delta[h^*\ot h]=[h^*_{(1)}\ot h_{(1)}]\ot [h^*_{(2)}\ot h_{(2)}]
\eqlabel{4.3.1}\\
&&\hspace*{-2cm} \varepsilon[h^*\ot h]=\lan h^*,\varepsilon_t(h)\ran\eqlabel{4.3.2}\\
&&\hspace*{-2cm} S[h^*\ot h]=[S^{-1}(h_{(2)}^*)\ot S(h_{(2)})]
\lan h_{(1)}^*,h_{(1)}\ran \lan h_{(3)}^*,S(h_{(3)})\ran \eqlabel{4.3.3}
\end{eqnarray}

\begin{proposition}\prlabel{4.3}
The $k$-linear isomorphism
$$f:\ H\bowtie H^*\to H^*\ot H,~~f(h\bowtie h^*)=h^*\ot S^{-1}(h)$$
is anti-multiplicative, and induces an algebra isomorphism
$f:\ D(H)\to D'(H)^{\rm op}$.
\end{proposition}

\begin{proof}
Let us first prove that $f$ reverses the multiplication. Indeed,
\begin{eqnarray*}
&&\hspace*{-2cm}
f(k\bowtie k^*)f(h\bowtie h^*)= (k^*\ot S^{-1}(k))(h^*\ot S^{-1}(h))\\
&=& (S^{-1}(k_{(1)})\leftact h^*\rightact k_{(3)})*k^*\ot S^{-1}(k_{(2)})S^{-1}(h)\\
&=& f((h\bowtie h^*)(k\bowtie k^*)).
\end{eqnarray*}
Using \leref{4.1a}, we easily compute that $f(J)=I$, and the result follows.
\end{proof}

Let us now define a comultiplication, counit and antipode on $D(H)$, in
such a way that $f:\ D(H)\to D'(H)$ is an isomorphism of Hopf algebras.
Obviously, the comultiplication is given by the formula
\begin{equation}\eqlabel{4.4.1}
\Delta[h\bowtie h^*]=[h_{(2)}\bowtie h^*_{(1)}]\ot [h_{(1)}\bowtie h^*_{(2)}].
\end{equation}
The counit is computed as follows:
\begin{equation}
\varepsilon[h\bowtie h^*]=\varepsilon[h^*\ot S^{-1}(h)]
\equal{\equref{4.3.2}}\lan h^*,\varepsilon_t(S^{-1}(h))\ran
\equal{\equref{1.4.1}}\lan h^*,1_{(2)}\ran \lan \varepsilon,h1_{(1)}\ran.
\eqlabel{4.4.2}
\end{equation}
Since the antipode of $H$ is the inverse of the antipode of $H^{\rm op}$,
the antipode of $D'(H)$ is transported to the inverse of the antipode of
$D(H)$. We find
\begin{eqnarray}
&&\hspace*{-2cm}S^{-1}[h\bowtie h^*]=
(f^{-1}\circ S\circ f)[h\bowtie h^*]=
f^{-1}(S[h^*\ot S^{-1}(h)])\nonumber\\
&=& f^{-1}[S^{-1}(h_{(2)}^*)\ot h_{(2)}]\lan h_{(1)}^*,S^{-1}(h_{(3)})\ran
\lan h_{(3)}^*,h_{(1)}\ran\nonumber\\
&=& [S(h_{(2)})\bowtie S^{-1}(h_{(2)}^*)]\lan h_{(1)}^*,S^{-1}(h_{(3)})\ran
\lan h_{(3)}^*,h_{(1)}\ran\eqlabel{4.4.3}
\end{eqnarray}
The antipode $S$ is then given by the formula
\begin{equation}\eqlabel{4.4.4}
S[h\bowtie h^*]=[S^{-1}(h_{(2)})\bowtie S(h_{(2)}^*)]
\lan h_{(1)}^*,S^{-1}(h_{(3)})\ran \lan h_{(3)}^*,h_{(1)}\ran
\end{equation}
Indeed,
\begin{eqnarray*}
&&\hspace*{-15mm}
S(S^{-1}[h\bowtie h^*])\\
&=& [h_{(3)}\bowtie h_{(3)}^*]
\lan h_{(1)}^*,S^{-1}(h_{(5)})\ran \lan h_{(2)}^*,h_{(4)}\ran
\lan h_{(5)}^*,h_{(1)}\ran \lan h_{(4)}^*,S^{-1}(h_{(2)})\ran\\
&=& [h_{(3)}\bowtie h_{(2)}^*]
\lan h_{(1)}^*,S^{-1}(h_{(5)})h_{(4)}\ran
\lan h_{(2)}^*,S^{-1}(h_{(2)})h_{(1)}\ran\\
&=& [h_{(2)}\bowtie h_{(2)}^*]
\lan h_{(1)}^*,\varepsilon_t(S^{-1}(h_{(3)}))\ran
\lan h_{(3)}^*,\varepsilon_t(S^{-1}(h_{(1)}))\ran\\
&=&\varepsilon([h\bowtie h^*]_{(1)})[h\bowtie h^*]_{(2)}\varepsilon([h\bowtie h^*]_{(3)})
= [h\bowtie h^*].
\end{eqnarray*}
Similar arguments show that $S^{-1}(S[h\bowtie h^*])=[h\bowtie h^*]$.

\begin{proposition}\prlabel{4.4}
Let $H$ be a weak Hopf algebra with bijective antipode, which is finitely
generated and projective as a $k$-module. Then
$D(H)$ is a weak Hopf algebra, with comultiplication, counit and antipode given by
the formulas (\ref{eq:4.4.1},\ref{eq:4.4.2},\ref{eq:4.4.3}). As a weak Hopf algebra,
$D(H)$ is isomorphic to $D'(H)^{\rm op}$.
\end{proposition}

\begin{proposition}\prlabel{4.5}
Let $H$ be a weak Hopf algebra with bijective antipode, which is finitely
generated and projective as a $k$-module.
The functor 
$$F:\ {}_H\YD^H\to {}_{D(H)}\ol{\Mm},~~F(M)=M,$$
with
$$(h\bowtie h^*)m=\lan h^*,m_{[1]}\ran hm_{[0]},$$
for all $h\in H$, $h^*\in H^*$ and $m\in M$ is an isomorphism of monoidal
categories.
\end{proposition}

\begin{proof}
We  already know (see \equref{1.42}) that $F$ is an isomorphism of categories,
so we only have to show that $F$ preserves the product. Take $M,N\in {}_H\YD^H$.
The right $H$-coaction on $M\ot_t N$ is given by the formula
(use \equref{2.5.1} and \equref{2.6.0}):
$$\rho(1_{(1)}m\ot 1_{(2)}n)=
m_{[0]}\ot n_{[0]}\ot n_{[1]}m_{[1]},$$
hence the left $D(H)$-action on $F(M\ot_t N)$ is the following
\begin{equation}\eqlabel{4.5.1}
[h\bowtie h^*](1_{(1)}m\ot 1_{(2)}n)= \lan h^*,n_{[1]}m_{[1]}\ran
h_{(1)}m_{[0]}\ot h_{(2)}n_{[0]}.
\end{equation}
We now compute 
$$F(N)\ot_t F(M)=\{[1\bowtie \varepsilon]X~|~X\in F(N)\ot F(M)\}.$$
Observe that
\begin{eqnarray*}
&&\hspace*{-2cm}
[1\bowtie\varepsilon]_{(1)}n\ot [1\bowtie\varepsilon]_{(2)}m
=\lan  \varepsilon_{(1)}, n_{[1]}\ran 1_{(2)}n_{[0]}\ot 
\lan \varepsilon_{(2)}, m_{[1]}\ran 1_{(1)}m_{[0]}\\
&=&\lan\varepsilon,n_{[1]}m_{[1]}\ran 1_{(2)}n_{[0]}\ot 1_{(1)}m_{[0]}.
\end{eqnarray*}
We claim that the switch map $\tau:\ M\ot N\to N\ot M$ induces an isomorphism
$\tau:\ F(M\ot_t N)\to F(N)\ot_t F(M)$ of $k$-modules. Indeed, take
$1_{(1)}m\ot 1_{(2)}n\in M\ot_t N$. Since $M\ot_t N$ is a Yetter-Drinfeld
module, we have that
$\varepsilon(n_{[1]}m_{[1]})m_{[0]}\ot n_{[0]}=1_{(1)}m\ot 1_{(2)}n$,
and
\begin{eqnarray*}
&&\hspace*{-2cm}
\tau(1_{(1)}m\ot 1_{(2)}n)
=1_{(2)}n\ot 1_{(1)}m= 1_{(2')}1_{(2)}n\ot 1_{(1')}1_{(1)}m\\
&=& \varepsilon(n_{[1]}m_{[1]})1_{(2)} n_{[0]}\ot 1_{(1)}m_{[0]}\\
&=&[1\bowtie\varepsilon]_{(1)}n\ot [1\bowtie\varepsilon]_{(2)}m
\in F(N)\ot_t F(M).
\end{eqnarray*}
Conversely,
$$\tau([1\bowtie\varepsilon]_{(1)}n\ot [1\bowtie\varepsilon]_{(2)}m)=
\varepsilon(n_{[1]}m_{[1]})1_{(1)}m_{[0]}\ot 1_{(2)} n_{[0]}
\in F(M\ot_t N).$$
Let us now show
that $\tau$ is left $D(H)$-linear. To this end, we compute the left
$D(H)$-action on $F(N)\ot_t F(M)$.
\begin{eqnarray*}
&&\hspace*{-15mm}
[h\bowtie h^*]\tau(1_{(1)}m\ot 1_{(2)}n)=
[h\bowtie h^*](1_{(2)}n\ot 1_{(1)}m)\\
&=&[h_{(2)}\bowtie h_{(1)}^*](1_{(2)}n)\ot [h_{(1)}\bowtie h_{(2)}^*](1_{(1)}m)\\
&\equal{\equref{2.7a.1}}&
\lan h_{(1)}^*,n_{[1]}S^{-1}(1_{(2)})\ran h_{(2)}n_{[0]}\ot
\lan h_{(2)}^*,1_{(1)}m_{[1]}\ran h_{(1)}m_{[0]}\\
&=& \lan h^*, n_{[1]}S^{-1}(1_{(2)})1_{(1)}m_{[1]}\ran 
h_{(2)}n_{[0]}\ot h_{(1)}m_{[0]}\\
&\equal{ \equref{4.5.1}}& \tau\bigl([h\bowtie h^*](1_{(1)}m\ot 1_{(2)}n)\bigr)
\end{eqnarray*}
It also follows that $F(H_t)$ is a unit object in ${}_{D(H)}\ol{\Mm}$.
Since the unit object in a monoidal category is unique up to automorphism,
we concluce that the target space of $D(H)_t$ is isomorphic to $H_t$.
This can also be seen as follows: in \cite{NVT}, it is shown that
$D'(H)_t=[\varepsilon\ot H_t]\cong H_t$. Since the target spaces of a weak
Hopf algebra and its opposite coincide,
it follows that $D(H)_t\cong H_t$.
\end{proof}

\section{Duality}\selabel{5}
Let $H$ be a weak Hopf algebra with bijective antipode, and
${}_H{\rm Rep}$ the category of left $H$-modules $M$ which are finitely
generated projective as a $k$-module. Let $M\in {}_H{\rm Rep}$, and
let $\{(n_i,n_i^*)~|~i=1,\cdots n\}$ be a finite dual basis of $M$.
From \cite{NVT}, we recall the following result. We refer to \cite{Kassel}
for the definition of duality in a monoidal category.

\begin{proposition}\prlabel{5.1}
The category ${}_H{\rm Rep}$ has left duality. The left dual of
$M\in {}_H{\rm Rep}$ is $M^*=\Hom(M,k)$ with left $H$-action defined by
\begin{equation}\eqlabel{5.1.1}
\lan h\cdot m^*, m\ran = \lan m^*,S(h)m\ran,
\end{equation}
for all $h\in H$, $m\in M$ and $m^*\in M^*$. The evaluation
map ${\rm ev}_M:\ M^*\ot_t M\to H_t$ and the coevaluation
map ${\rm coev}_M:\ H_t\to M\ot_t M^*$ are defined as follows:
$${\rm ev}_M(1_{(1)}\cdot m^*\ot 1_{(2)}m)=\lan m^*,1_{(1)}m\ran 1_{(2)};$$
$${\rm coev}_M(z)=z\cdot (\hbox{$\sum_i$} n_i\ot n_i^*).$$
\end{proposition}

Let $M$ be a finitely generated projective left $H$-comodule. Then $M^*$
is also a left $H$-comodule, with left $H$-coaction $\lambda:\ M^*\to
H\ot M^*$ given by
$$\lambda(m^*)=\hbox{$\sum_i$} \lan m^*,n_{i[0]}\ran S^{-1}(n_{i[-1]})\ot n_i^*.$$
The definition of $\lambda$ can also be stated as follows:
$\lambda(m^*)=m^*_{[-1]}\ot m^*_{[0]}$ if and only if
\begin{equation}\eqlabel{5.2.1}
\lan m_{[0]}^*,m\ran S(m_{[-1]}^*)=\lan m^*,m_{[0]}\ran m_{[-1]},
\end{equation}
for all $m\in M$.

\begin{proposition}\prlabel{5.2}
Let $M$ be a finitely generated projective left-left Yetter-Drinfeld module
over the weak Hopf algebra $H$. Then $M^*$ with $H$-action and $H$-coaction
given by \equref{5.1.1} and \equref{5.2.1} is also a left-left Yetter-Drinfeld
module.
\end{proposition}

\begin{proof}
We have to show that
$$\lambda(h\cdot m^*)=\hbox{$\sum_i$} \lan m^*,S(h)n_{i[0]}\ran S^{-1}(n_{i[-1]})\ot n_i^*$$
equals
$$h_{(1)}m_{[-1]}^*S(h_{(3)})\ot h_{[2]}m_{[-1]}^*\\
= \hbox{$\sum_i$} \lan m^*,n_{i[0]}\ran h_{(1)}S^{-1}(n_{i[-1]})S(h_{(3)})\ot
(h_{(2)}\cdot n_i^*).$$
It suffices to show that both terms coincide after we evaluate the second tensor
factor at an arbitrary $m\in M$.
\begin{eqnarray*}
&&\hspace*{-2cm}
\hbox{$\sum_i$} \lan m^*,n_{i[0]}\ran h_{(1)}S^{-1}(n_{i[-1]})S(h_{(3)})
\lan n_i^*,S(h_{(2)})m\ran\\
&=&\lan m^*,(S(h_{(2)})m)_{[0]}\ran h_{(1)}
S^{-1}\Bigl((S(h_{(2)})m)_{[-1]}\Bigr)S(h_{(3)})\\
&\equal{\equref{2.2.1}}&
\lan m^*, S(h_{(3)})m_{[0]}\ran h_{(1)}
S^{-1}\Bigl(S(h_{(4)})m_{[-1]}S^2(h_{(2)})\Bigr) S(h_{(5)})\\
&=& \lan m^*, S(h_{(3)})m_{[0]}\ran h_{(1)}
S(h_{(2)})S^{-1}(m_{[-1]})h_{(4)} S(h_{(5)})\\
&=& \lan m^*, S(h_{(2)})m_{[0]}\ran \varepsilon_t(h_{(1)})
S^{-1}(m_{[-1]})\varepsilon_t(h_{(3)})
\end{eqnarray*}
\begin{eqnarray*}
&\equal{\equref{1.6.3}}&
\lan m^*, S(1_{(2)}h_{(1)})m_{[0]}\ran S(1_{(1)})
S^{-1}(m_{[-1]})\varepsilon_t(h_{(2)})\\
&=& \lan m^*, S(h_{(1)})1_{(1)}m_{[0]}\ran 1_{(2)}
S^{-1}(m_{[-1]})\varepsilon_t(h_{(2)})\\
&\equal{(\ref{eq:1.1.3},\ref{eq:2.3.1})}&
\lan m^*, S(1_{(1)}h)m_{[0]}\ran S^{-1}(m_{[-1]})1_{(2)}\\
&=& \lan m^*, S(h)S(1_{(1)})m_{[0]}\ran S^{-1}(S(1_{(2)})m_{[-1]})\\
&=& \lan m^*, S(h)1_{(2)}m_{[0]}\ran S^{-1}(1_{(1)}m_{[-1]})\\
&\equal{\equref{2.1.1}}&
\lan m^*, S(h)m_{[0]}\ran S^{-1}(m_{[-1]})\\
&=& \hbox{$\sum_i$}\lan m^*,S(h)n_{i[0]}\ran S^{-1}(n_{i[-1]}^*)\lan n_i^*,m\ran
\end{eqnarray*}
\end{proof}

\begin{proposition}\prlabel{5.3}
The category of finitely generated projective left-left Yetter-Drinfeld modules
has left duality.
\end{proposition}

\begin{proof}
In view of the previous results, it suffices to show that the evaluation map
${\rm ev}_M$ and the coevaluation map ${\rm coev}_M$ are left $H$-colinear,
for every finitely generated projective left-left Yetter-Drinfeld module $M$.
Let us first show that ${\rm ev}_M$ is left $H$-colinear.
\begin{eqnarray*}
&&\hspace*{-2cm}
(H\ot {\rm ev}_M)(\lambda(1_{(1)}\cdot m^*\ot 1_{(2)}m))\\
&=& m^*_{[-1]}m_{[-1]}\ot \lan m^*_{[0]},1_{(1)}m_{[0]}\ran 1_{(2)}\\
&\equal{\equref{5.2.1}}&
\lan m^*,(1_{(1)}m_{[0]})_{[0]}\ran S^{-1}((1_{(1)}m_{[0]})_{[-1]})m_{[-1]}\ot 1_{(2)}\\
&\equal{\equref{2.2bis.1}}&
\lan m^*,m_{[0]}\ran 1_{(1)}S^{-1}(m_{[-1]})m_{[-2]}\ot 1_{(2)}\\
&=&\lan m^*,m_{[0]}\ran 1_{(1)}\varepsilon_t(S^{-1}(m_{[-1]}))\ot 1_{(2)}\\
&\equal{\equref{2.3.1}}&
\lan m^*,1_{(1')}m_{[0]}\ran 1_{(1)}\varepsilon_t(1_{(2')}S^{-1}(m_{[-1]}))\ot 1_{(2)}\\
&\equal{\equref{1.1.7}}&
\lan m^*,1_{(1')}m_{[0]}\ran 1_{(1)}1_{(2')}\varepsilon_t(S^{-1}(m_{[-1]}))\ot 1_{(2)}\\
&\equal{\equref{1.9.2}}&
\lan m^*,1_{(1')}S^{-1}(\varepsilon_t(S^{-1}(m_{[-1]})))m_{[0]}\ran 1_{(1)}1_{(2')}\ot 1_{(2)}\\
&\equal{\equref{1.7.1}}&
\lan m^*,1_{(1')}\varepsilon_s(S^{-2}(m_{[-1]}))m_{[0]}\ran 1_{(1)}1_{(2')}\ot 1_{(2)}\\
&\equal{(\ref{eq:1.1.1},\ref{eq:2.3.2})}&
\lan m^*,1_{(1)}m\ran 1_{(2)}\ot 1_{(3)}\equal{\equref{2.5.2}}
\lambda(\lan m^*,1_{(1)}m\ran 1_{(2)})\\
&=&\lambda({\rm ev}_M(1_{(1)}\cdot m^*\ot 1_{(2)}m)).
\end{eqnarray*}
To prove that ${\rm coev}_M$ is left $H$-colinear, we have to show that,
for all $z\in H_t$,
\begin{eqnarray*}
&&\hspace*{-2cm}
\lambda({\rm coev}_M(z))=\hbox{$\sum_i$} \lambda(1_{(1)}zn_i\ot 1_{(2)}\cdot n_i^*)\\
&=&
\hbox{$\sum_i$} (1_{(1)}zn_i)_{[-1]}(1_{(2)}\cdot n_i^*)_{[-1]}
\ot (1_{(1)}zn_i)_{[0]}\ot (1_{(2)}\cdot n_i^*)_{[0]}
\end{eqnarray*}
equals
$$
(H\ot {\rm coev}_M)(\lambda(z))=
(H\ot {\rm coev}_M)(1_{(1)}z\ot 1_{(2)})=
 \hbox{$\sum_i$} 1_{(1)}z\ot 1_{(2)}n_i\ot 1_{(3)}\cdot n_i^*.
$$
It suffices to show that both terms coincide after we evaluate the third tensor
factor at an arbitrary $m\in M$. Indeed
\begin{eqnarray*}
&&\hspace*{-2cm}
\hbox{$\sum_i$} (1_{(1)}zn_i)_{[-1]}(1_{(2)}\cdot n_i^*)_{[-1]}
\ot (1_{(1)}zn_i)_{[0]}\lan (1_{(2)}\cdot n_i^*)_{[0]},m\ran\\
&\equal{\equref{5.2.1}}&
\hbox{$\sum_i$} (1_{(1)}zn_i)_{[-1]} \lan 1_{(2)}\cdot n_i^*, m_{[0]}\ran S^{-1}(m_{[-1]})
\ot (1_{(1)}zn_i)_{[0]}\\
&=&
\hbox{$\sum_i$} (1_{(1)}zn_i)_{[-1]} \lan n_i^*, S(1_{(2)})m_{[0]}\ran S^{-1}(m_{[-1]})
\ot (1_{(1)}zn_i)_{[0]}\\
&=& (1_{(1)}zS(1_{(2)})m_{[0]})_{[-1]} S^{-1}(m_{[-1]}) \ot (1_{(1)}zS(1_{(2)})m_{[0]})_{[0]} 
\end{eqnarray*}
\begin{eqnarray*}
&\equal{(\ref{eq:1.1.4},\ref{eq:2.2.1})}&
1_{(1)}zm_{[-1]}S(1_{(3)})S^{-1}(m_{[-2]}) \ot 1_{(2)}m_{[0]}\\
&=&
1_{(1)}zm_{[-1]}S(1_{(2')})S^{-1}(m_{[-2]}) \ot 1_{(2)}1_{(1')}m_{[0]}\\
&\equal{\equref{2.2.2}}&
1_{(1)}zm_{[-1]}S^{-1}(m_{[-2]}) \ot 1_{(2)}m_{[0]}\\
&=&
1_{(1)}zS^{-1}(\varepsilon_t(m_{[-1]})) \ot 1_{(2)}m_{[0]}\\
&\equal{\equref{1.9.2}}&
1_{(1)}z \ot 1_{(2)}\varepsilon_t(m_{[-1]})m_{[0]}
\equal{\equref{2.1.3}}1_{(1)}z \ot 1_{(2)}m\\
&=& 1_{(1)}z \ot 1_{(2)}S(1_{(3)})m=
\hbox{$\sum_i$} 1_{(1)}z \ot 1_{(2)}n_i\lan n_i^*,S(1_{(3)})m\ran\\
&=& \hbox{$\sum_i$} 1_{(1)}z \ot 1_{(2)}n_i\lan 1_{(3)}\cdot n_i^*,m\ran,
\end{eqnarray*}
as needed.
\end{proof}

\section{Appendix. Weak bialgebras and bialgebroids}\selabel{6}
In \cite{Schauenburg0}, Yetter-Drinfeld modules over a $\times_R$-bialgebra
(see \cite{Takeuchi}) are introduced,
and it is shown that the weak center of the category of left modules 
is isomorphic to the category of Yetter-Drinfeld modules. The notion of $\times_R$-bialgebra
is equivalent to the notion of $R$-bialgebroid, we refer to \cite{BM} for a detailed
discussion. So we can consider Yetter-Drinfeld modules over bialgebroids.\\
A weak bialgebra $H$ can be viewed as a bialgebroid over the target space
$H_t$; this was shown in \cite{EN} in the weak Hopf algebra case, and generalized
to the weak bialgebra case in \cite{Schauenburg}. The aim of this Section is to make
clear that Yetter-Drinfeld modules over $H$ considered as a weak bialgebra
coincide with Yetter-Drinfeld modules over $H$-considered as a bialgebroid.\\
To this end, we first recall the definition of a bialgebroid, as introduced by Lu
\cite{Lu}. Let $k$ be a commutative ring, and $R$ a $k$-algebra. An $R\ot R^{\rm op}$-ring
is a pair $(H,i)$, with $H$ a $k$-algebra and $i:\ R\ot R^{\rm op}\to H$.
Giving $i$ is equivalent to giving algebra maps $s_H:\ R\to H$ and
$t_H:\ R\to H^{\rm op}$ satisfying $s_H(a)t_H(b)=t_H(b)s_H(a)$, for all $a,b\in R$.
We then have that $i(a\ot b)=s_H(a)t_H(b)$. Restriction of scalars makes $H$
into a left $R\ot R^{\rm op}$-module, and an $R$-bimodule:
$$a\cdot h\cdot b=s_H(a)t_H(b)h.$$
Consider
\begin{eqnarray*}
&&\hspace*{-2cm}
H\times_R H=\{\sum_i h_i\ot_R k_i\in H\ot_R H\\
&|&
\sum_i h_it_H(a)\ot_R k_i=\sum_i h_i\ot_R k_is_H(a),~~{\rm for~all~}a\in R\}
\end{eqnarray*}
It is easy to show that $H\times_R H$ is a $k$-subalgebra of $H\ot_R H$.\\
Recall that an $R$-coring is a triple $(H,\tilde{\Delta},\tilde{\varepsilon})$,
with $H$ an $R$-bimodule and $\tilde{\Delta}:\ H\to H\ot_R H$ and
$\tilde{\varepsilon}:\ H\to R$ $R$-bimodule maps satisfying the usual coassociativity
and counit properties; we refer to \cite{BW} for a detailed discussion of corings.

\begin{definition}\delabel{6.1} \cite{Lu}
A left $R$-bialgebroid is a fivetuple $(H,s_H,t_H,\tilde{\Delta},\tilde{\varepsilon})$
satisfying the following conditions.
\begin{enumerate}
\item $(H,\tilde{\Delta},\tilde{\varepsilon})$ is an $R$-coring;
\item $(H,m\circ (s_H\ot t_H)=i)$ is an $R\ot R^{\rm op}$-ring;
\item $\im(\tilde{\Delta})\subset H\times_R H$;
\item $\tilde{\Delta}:\ H\to H\times_R H$ is an algebra map, $\tilde{\varepsilon}(1_H)=
1_R$ and
$$\tilde{\varepsilon}(gh)=\tilde{\varepsilon}(gs_H(\tilde{\varepsilon}(h)))=
\tilde{\varepsilon}(gt_H(\tilde{\varepsilon}(h))),$$
for all $g,h\in H$.
\end{enumerate}
\end{definition}

Take two left $H$-modules $M$ and $N$; then $M$ and $N$ are $R$-bimodules,
by restriction of scalars. $M\ot_RN$ is a left $H$-module, with
$$h\cdot (m\ot_Rn)=h_{(1)}m\ot_R h_{(2)}n.$$
Also $R$ is a left $H$-module, with
$$h\cdot r=\tilde{\varepsilon}(hs_H(r))=\tilde{\varepsilon}(ht_H(r)).$$
$({}_H\Mm,\ot_R,R)$ is a monoidal category, and the restriction of scalars
functor ${}_H\Mm\to {}_R\Mm_R$ is strictly monoidal; this can be used to reformulate
the definition of a bialgebroid (see \cite{BM,Schauenburg98,Szlachanyi}).\\
In \cite[Sec. 4]{Schauenburg0}, left-left Yetter-Drinfeld modules over $H$ are introduced,
and it is shown that $\Ww_l({}_H\Mm)$ is isomorphic to the category of Yetter-Drinfeld
modules. According to \cite{Schauenburg0}, a left-left Yetter-Drinfeld $H$-module
is a left comodule $M$ over the coring $H$, together with a left $H$-action on $M$
such that the underlying left $R$-actions coincide, and such that 
\begin{equation}\eqlabel{6.1.1}
h_{(1)}m_{[-1]}\ot_R h_{(2)}\cdot m_{[0]}=(h_{(1)}\cdot m)_{[-1]}h_{(2)}
\ot_R (h_{(1)}\cdot m)_{[0]}
\end{equation}
holds in $H\ot_R M$, for all $h\in H$ and $m\in M$.\\

Let $H$ be a weak bialgebra, and consider the maps
\begin{eqnarray*}
&&\hspace*{-2cm} s_H:\ H_t\rTo^{\subset} H;\\
&&\hspace*{-2cm} t_H=\ol{\varepsilon}_{s|H_t}:\ H_t\to H_s\subset H;\\
&&\hspace*{-2cm} \tilde{\Delta}={\rm can}\circ \Delta:\
H\to H\ot H\rTo^{\rm can}H\ot_{H_t}H;\\
&&\hspace*{-2cm} \tilde{\varepsilon}=\varepsilon_t:\ H\to H_t.
\end{eqnarray*}
Then $(H,s_H,t_H,\tilde{\Delta},\tilde{\varepsilon})$ is a left $H_t$-bialgebroid.
The fact that $\im(\tilde{\Delta})\subset H\times_{H_t}H$ follows from the separability
of $H_t$ as a $k$-algebra (cf. \prref{1.2b}).\\
We have seen in \seref{1.1} that, for any two left $H$-modules $M$ and $N$,
we have an isomorphism $\ol{\pi}:\ M\ot_{H_t}N\to M\ot_t N$. This entails that
the monoidal categories $({}_H\Mm,\ot_t,H_t)$ and $({}_H\Mm,\ot_{H_t},H_t)$ are
isomorphic, and a fortiori, their weak left centers are isomorphic categories.
Consequently, the two corresponding categories of Yetter-Drinfeld modules are
isomorphic. This can also be seen directly, comparing the definitions in \seref{2}
and \equref{6.1.1}.

\begin{center}
{\bf Acknowledgment}
\end{center}
We thank Tomasz Brzezi\'nski and the referee for their useful comments, and Adriana Nenciu
for sending us her paper \cite{Nenciu}.

\end{document}